\documentclass[12pt]{article}
\usepackage{latexsym}
\usepackage{epsfig}
\usepackage{amsfonts,amssymb,amsmath}
\usepackage{stmaryrd}
\usepackage[latin1]{inputenc}
\usepackage[T1]{fontenc}

\title{Incremental moments and H\"older exponents \\
of multifractional multistable processes}
\author{R. Le Gu\'evel\\
\small{{\em Universit\'{e} Nantes, Laboratoire de Math\'{e}matiques 
Jean Leray UMR CNRS 6629}}\\
\small{{\em 2 Rue de la Houssini\`{e}re - BP 92208 - F-44322 Nantes Cedex 3, France}}\\
\small{{\em ronan.leguevel@univ-nantes.fr}}\\
\small{ and } \\
J. L\'{e}vy V\'{e}hel \\
\small{{\em Regularity team, INRIA Saclay, Parc Orsay Universit\'e}}\\
\small{{\em 4 rue Jacques Monod - Bat P - 91893 Orsay Cedex, France}}\\
\small{{\em jacques.levy-vehel@inria.fr}}}
\date{}
\def\bbbr{{\bf R}} 

\newtheorem{theo}{Theorem}

\newtheorem{prop}[theo]{Proposition}

\newtheorem{lem}[theo]{Lemma}

\newcommand\E{\mbox{\sf E}}

\newcommand\efd{\stackrel{fdd}{=}}

\newcommand\ed{\stackrel{d}{=}}
\renewcommand\P{{\sf P}}

\newcommand\h{{\mathcal H}_t}

\newcommand{\one}{\ifmmode {\sf 1}\hspace{-.26em}{\sf
l}\hspace{-.35em}{\sf \_} \else ${\sf 1}\hspace{-.26em}{\sf
l}\hspace{-.35em}{\sf \_}$ \fi}

\newcommand\varep{\varepsilon}

\newcommand\ho{H\"{o}lder }

\renewcommand{\Box}{\mbox{\rule{1ex}{1ex}}}

\begin{document}

\maketitle

\begin{abstract}
\noindent Multistable processes, that is, processes which are, at each
``time'', tangent to a stable process, but where the index of stability
varies along the path, have been recently introduced as models for
phenomena where the intensity of jumps is non constant. In this work,
we give further results on (multifractional) multistable processes related to their local
structure. We show that, under certain conditions, the incremental
moments display a scaling behaviour, and that the pointwise \ho
exponent is, as expected, lower than the localisability index.
    
\end{abstract}

\vspace{1cm}


\vspace{1cm}


\vspace{1cm}
\section{Introduction}

Multistable processes are stochastic processes which are ``locally stable'', but where
the index of stability varies with ``time''. To be more precise, we need to recall the
definition of a {\it localisable process} \cite{Fal5,Fal6}: $Y =\{Y(t): t \in \bbbr\}$ 
is said to be $h-$localisable at $u$ if there exists an $h \in \bbbr$ 
and a non-trivial limiting process $Y_{u}'$ such that
\begin{equation}
\lim_{r \to 0}\frac{Y(u+rt) -Y(u)}{r^{h}} = Y_{u}'(t).
\label{locform1}
\end{equation}
(Note $Y_{u}'$ may and in general will vary with $u$.) When the limit exits, 
$Y_{u}'=\{Y_{u}'(t): t\in \bbbr\}$ is termed the {\it local form} or tangent
process of $Y$ at $u$.
The limit (\ref{locform1}) may
be taken in mainly two ways: convergence 
in finite dimensional distributions, or in distribution
(in which case the process is called
{\it strongly $h$-localisable}). In the sequel, equality in 
finite dimensional distributions will be denoted $\efd$, and
equality in distributions $\ed$.

A classical example of localisable process is 
multifractional Brownian motion $Y$ \cite{AL2,BJR,EH,PL}
which ``looks like'' index-$h(u)$ fractional
Brownian motion close to time $u$ but where $h(u)$ varies, that is  
\begin{equation}
\lim_{r \to 0}\frac{Y(u+rt) -Y(u)}{r^{h}} = B_{h(u)}(t) \label{exfbm}
\end{equation}
where $B_{h}$ is index-$h$ fractional Brownian motion.
A generalization of mBm, where the Gaussian
measure is replaced by an $\alpha$-stable one, leads to
multifractional stable processes, where the
local form is an $h(u)$-self-similar linear
$\alpha$-stable motion \cite{ST1,ST2}.

Multifractional multistable processes provide a further step of generalization: they are localisable processes such that the tangent
process is again an $\alpha$-stable random process, but where $\alpha$
now varies with time.
Multifractional multistable processes were constructed in 
\cite{FLGLV,FLV,LGLV,FL} using respectively 
moving averages, sums over Poisson processes, the Ferguson - Klass - LePage series 
representation, and multistable measures. Section \ref{Ex} below provides several 
specific examples of such processes.

\bigskip

The aim of this work is twofold:
\begin{enumerate}
	\item We show that, for a large class of (multifractional) multistable 
processes, a precise estimate for the incremental moments holds. More precisely, we 
prove in section \ref{moments} that there exists a natural scaling relation for $\E \left[ |Y(t+\varep)-Y(t)|^{\eta}\right]$ and $\varep$ small. This class includes (multifractional) multistable processes considered in \cite{FLV,LGLV}, in particular Lévy
multistable motions and linear multistable multifractional motions.

	\item We then study the pointwise \ho regularity of (multifractional) multistable 
	processes. For the same
	class as above, we obtain an almost sure upper bound for this exponent. 
	In the case of the Lévy multistable motion, we are able to compute its exact value. Not surprisingly, it turns out to be equal, at each point, almost surely, to the localisability index. 
	Note that a uniform statement, {\it i.e.} a statement like ``almost surely, at
	each point'', cannot hold true in general. Indeed, it already fails for the 
	case of a Lévy stable motion. The right frame in this respect is multifractal analysis,
	and results in this direction will be presented in a forthcoming work.
\end{enumerate}

\smallskip

The remainder of this work is organized as follows. In the next section, we recall the definition
of multistable processes based on the Ferguson - Klass - LePage series representation 
used in \cite{LGLV} (this defines processes
which are equal in distribution to the ones obtained in \cite{FLV} through sums over Poisson 
processes). Our main results on incremental moments and upper bound for the pointwise Hölder exponents are described in section \ref{MR}. Subsection \ref{Ex} applies these results to the
linear multistable multifractional motion. In
subsection \ref{exlmm}, we state the fact that, for the Lévy multistable motion, the exact almost
sure value of the exponent is indeed the reciprocal of the localisability exponent.
In section \ref{IR}, we give intermediate
results, some of which being of independent interest, which are used in the proofs of the main
statements. Section \ref{proo} gathers technical results followed by the proofs of the statements
related with the incremental moments and upper bounds on the exponents. Section \ref{proolevy}
contains the proof of the lower bound for the exponent of the multistable Lévy motion.
Finally, section \ref{Assum} gives a list of the various technical conditions on multistable processes required  by our approach so that their incremental moments and Hölder exponents may be estimated.

\section{Multistable processes}\label{MSP}

Our results will apply to certain processes that are defined as ``diagonals''
of random fields that we describe in subsections \ref{fmsc} and \ref{sfmsc}.

\subsection{Finite measure space case}\label{fmsc}

Let $(E,{\cal E},m)$ be a finite measure space, and $U$ be an open interval of $\bbbr$. Let $\alpha$ be a $C^1$ function defined on $U$ and ranging in $[c,d] \subset (0,2)$. Let $b$ be a $C^1$ function defined and bounded on $U$. Let $f(t,u,.)$ be a family of functions such that, for all $(t,u) \in U^2$, $f(t,u,.) \in {\cal F}_{\alpha(u)}(E,{\cal E}, m)$. Let $(\Gamma_i)_{i \geq 1}$ be a sequence of arrival times of a Poisson process with unit arrival time, $(V_i)_{i \geq 1}$ be a sequence of i.i.d. random variables with distribution $\hat m = m/m(E)$ on $E$, and $(\gamma_i)_{i \geq 1}$ be a sequence of i.i.d. random variables with distribution $P(\gamma_i=1)=P(\gamma_i=-1)=1/2$. Assume finally that the three sequences $(\Gamma_i)_{i \geq 1}$, $(V_i)_{i \geq 1}$, and $(\gamma_i)_{i \geq 1}$ are independent. As in \cite{LGLV}, we will
consider the following random field:

\begin{equation}\label{msfm}
X(t,u)= b(u) (m(E))^{1/\alpha(u)}C^{1/\alpha(u)}_{\alpha(u)} \sum_{i=1}^{\infty} \gamma_i \Gamma_i^{-1/\alpha(u)} f(t,u,V_i),
\end{equation}
where $C_{\eta} = \left( \int_{0}^{\infty} x^{-\eta} \sin (x)dx \right)^{-1}$.

Note that when the function $\alpha$  is constant, then \eqref{msfm}
is just the Ferguson - Klass - LePage series representation of a stable
random variable (see \cite{BJP,FK,LP1,LP2,JR}
and \cite[Theorem 3.10.1]{ST} for the specific properties of this representation that
will be needed here).

\subsection{$\sigma$-finite measure space case}\label{sfmsc}

When the space $E$ has infinite measure, one cannot use definition (\ref{msfm}), since it is no longer possible to renormalize by $m(E)$. However, in the $\sigma$-finite case, one may always perform a change of measure that allows to reduce to the finite case, as explained in \cite{ST} proposition 3.11.3 (see also section 4 of \cite{LGLV}). In our frame, this simply means adding a term
involving the change of measure in the definition of the field.

Let $(E,{\cal E},m)$ be a $\sigma$-finite measure space and $U$ be an open interval of $\bbbr$. Let $r : E \rightarrow \mathbb{R}_+$ be such that $\hat{m}(dx)=\frac{1}{r(x)}m(dx)$ is a probability measure. Let $\alpha$ be a $C^1$ function defined on $U$ and ranging in $[c,d] \subset (0,2)$. Let $b$ be a $C^1$ function defined and bounded on $U$. Let $f(t,u,.)$ be a family of functions such that, for all $(t,u) \in U^2$, $f(t,u,.) \in {\cal F}_{\alpha(u)}(E,{\cal E}, m)$. Let $(\Gamma_i)_{i \geq 1}$ be a sequence of arrival times of a Poisson process with unit arrival time, $(V_i)_{i \geq 1}$ be a sequence of i.i.d. random variables with distribution $\hat m $ on $E$, and $(\gamma_i)_{i \geq 1}$ be a sequence of i.i.d. random variables with distribution $P(\gamma_i=1)=P(\gamma_i=-1)=1/2$. Assume finally that the three sequences $(\Gamma_i)_{i \geq 1}$, $(V_i)_{i \geq 1}$, and $(\gamma_i)_{i \geq 1}$ are independent. We consider again a random field:

\begin{equation}\label{msfm2}
X(t,u)= b(u)C^{1/\alpha(u)}_{\alpha(u)} \sum_{i=1}^{\infty} \gamma_i \Gamma_i^{-1/\alpha(u)} r(V_i)^{1/\alpha(u)}f(t,u,V_i),
\end{equation}
with $C_{\alpha}$ as above.

\subsection{The diagonal processes}

Multistable processes are obtained by taking diagonals on $X$, {\it i.e.}
defining $Y(t) = X(t,t)$, both in the finite and $\sigma$-finite measure space cases. Indeed, as shown in Theorems 3.3 and 4.5 of \cite{LGLV}, provided some conditions are satisfied both by $X$ and by the function $f$, $Y$ will be
a localisable process whose local form is a stable process. In the remaining
of this work, we obtain, under certain assumptions (which imply that $Y$
is indeed localisable), estimates on the incremental
moments and the pointwise \ho regularity of $Y$.

\section{Main results}\label{MR}

The three following theorems apply to a diagonal process $Y$
defined from the field $X$ given by (\ref{msfm}) or (\ref{msfm2}). 
For convenience, the conditions required on $X$ and 
the function $f$ that appears in
(\ref{msfm}) and (\ref{msfm2}), denoted (C1), \ldots, (C15), are gathered in section \ref{Assum}.

\subsection{Moments of multistable processes}\label{moments}

\begin{theo}\label{vitesp}
 Let $t \in \bbbr$ and $U$ be an open interval of $\bbbr$ with $t \in U$. Let  $\eta \in (0,c)$. Suppose that $f$ satisfies (C1), (C2), (C3) (or (C1), (Cs2), (Cs3), (Cs4) in the $\sigma$-finite case), and (C9), and that $X$ verifies (C5) at $t$. Then, when $\varep$ tends to 0,
 \begin{displaymath}
 \E \left[ |Y(t+\varep)-Y(t)|^{\eta}\right] \sim \varep^{\eta h(t)} \E \left[ |Y'_t(1)|^{\eta}\right].
\end{displaymath}
\end{theo}

\noindent
{\bf Proof} 

\noindent
See section {\ref{proo}}.

\medskip

\noindent
{\bf Remark:} Under the conditions listed in the theorem, Theorems 3.3 and 
4.5 of \cite{LGLV} imply that $Y$ is $h(t)-$localisable at $t$.

\subsection{Pointwise H\"older exponent of multistable processes}\label{HE}

Let $\h = \sup \{ \gamma : \lim_{r \rightarrow 0}\limits \frac{| Y(t+r) - Y(t)|}{|r|^{\gamma}} = 0\}$ denote the \ho exponent of the (non-differentiable)
process $Y$ at $t$.





\begin{theo}(Upper bound)\label{upbound}
Suppose that there exists a function $h$ defined on $U$ such that (C6), (C7), (C8), (C10), (C11), (C12), (C13), (C14) and (C15) holds for some $t \in U$. Assuming (C1), (C2), (C3), (or (C1), (Cs2), (Cs3), (Cs4) in the $\sigma$-finite case), one has:
\begin{displaymath}
 \h \leq h(t).
\end{displaymath}
\end{theo}

\noindent
{\bf Proof} 

\noindent
See section {\ref{proo}}.

\subsection{Example: the linear multistable multifractional motion}\label{Ex}

In this section, we apply the results above to the ``multistable version'' of a 
classical process known as the linear stable multifractional motion, which is itself
a extension of the linear stable fractional motion, defined as follows (in the sequel, 
$M$ will always denote a symmetric $\alpha$-stable ($ 0 < \alpha < 2$) random measure on $\bbbr$ with control measure Lebesgue measure ${\cal L}$):

\begin{displaymath}
L_{\alpha,H,b^{+},b^{-}} (t) = \int_{-\infty}^{\infty} f_{\alpha,H}(b^{+},b^{-},t,x) M(dx) 
\end{displaymath}
where $t \in \mathbb{R}$, $H \in (0,1)$, $b^{+},b^{-} \in \mathbb{R}$, and
\begin{align*}
f_{\alpha,H}(b^{+},b^{-},t,x) = b^{+} & \left( (t-x)_{+}^{H - 1/ \alpha} - (-x)_{+} ^{H - 1/ \alpha} \right) \nonumber
\\
& + b^{-} \left( (t-x)_{-} ^{H- 1/ \alpha} - (-x)_{-} ^{H- 1/ \alpha} \right). \label{lfsm}
\end{align*}

When $b^{+}=b^{-}=1$, this process
is called well-balanced linear fractional $\alpha$-stable motion and denoted $L_{\alpha,H}$.

The localisability of the
linear fractional $\alpha$-stable motion simply stems from the fact that
it is $1/\alpha$-self-similar with stationary increments \cite{Fal6}. 

The multistable version of this processes was defined in \cite{FLGLV,FLV}. Its incremental
moments and regularity are described by the following theorems:

\begin{theo} {\it (Linear multistable multifractional motion)}.\label{esplmmm}
 Let   $\alpha:\bbbr \to [c,d] \subset (0,2)$ and
$H:\bbbr \to (0,1)$ be continuously differentiable. 
 Let $(\Gamma_i)_{i \geq 1}$ be a sequence of arrival times of a Poisson process with unit arrival time, 
$(V_i)_{i \geq 1}$ be a sequence of i.i.d. random variables with distribution $\hat{m}(dx) = \frac{3}{\pi^2}\sum_{j=1}^{+\infty} j^{-2} \mathbf{1}_{[-j,-j+1[ \cup [j-1,j[}(x) dx$ on $\bbbr$, and $(\gamma_i)_{i \geq 1}$ be a sequence of i.i.d. random variables with distribution $P(\gamma_i=1)=P(\gamma_i=-1)=1/2$. Assume finally that the three sequences $(\Gamma_i)_{i \geq 1}$, $(V_i)_{i \geq 1}$, and $(\gamma_i)_{i \geq 1}$ are independent and 
define
\begin{equation} 
X(t,u)=C_{\alpha(u)}^{1/ \alpha(u)} (\frac{\pi^2 j^2}{3})^{1/ \alpha(u)}\sum_{i,j=1}^{\infty} \gamma_i \Gamma_i^{-1/\alpha(u)} (|t-V_i|^{H(u)-1/ \alpha(u)}-|V_i|^{H(u)-1/\alpha(u)}) \mathbf{1}_{[-j,-j+1[ \cup [j-1,j[}(V_i)
\end{equation}
and the linear multistable multifractional motion
\begin{displaymath}
Y(t) = X(t,t).
\end{displaymath}
Then for all $t \in \bbbr$ and $\eta < c$, when $\varep$ tends to 0,
 \begin{displaymath}
 \E \left[ |Y(t+\varep)-Y(t)|^{\eta}\right] \sim \frac{2^{\eta-1} \Gamma(1-\frac{\eta}{\alpha(t)})}{\eta \int_{0}^{\infty} u^{-\eta-1} \sin^2(u)du} \left( \int_{\bbbr} \left| | 1-x |^{H(t)-\frac{1}{\alpha(t)}} - |x|^{H(t)-\frac{1}{\alpha(t)}}\right|^{\alpha(t)} dx \right)^{\frac{\eta}{\alpha(t)}} \varep^{\eta H(t)} .
\end{displaymath}
\end{theo}

\noindent
{\bf Proof} 

\noindent
See section {\ref{proo}}.

\begin{theo}\label{expolmmm}
Let $Y$ be the linear multistable multifractional motion defined on $\bbbr$ with $H-\frac{1}{\alpha}$ a non-negative function. For all $t \in \bbbr$, almost surely,
\begin{displaymath}
\h \leq H(t).
\end{displaymath}
\end{theo}

\noindent
{\bf Proof} 

\noindent
See section {\ref{proo}}.

\subsection{Example: the Lévy multistable motion}\label{exlmm}

In the case of the Lévy multistable motion, we are able to provide a 
more precise result, to the effect that, at each point,
the exact almost sure value of the Hölder exponent is known. Let us first recall some definitions. With $M$ again denoting a symmetric $\alpha$-stable ($ 0 < \alpha < 2$) random measure on $\bbbr$ with control measure Lebesgue measure ${\cal L}$, we write
$$L_{\alpha} (t) := \int_{0}^{t} M(dz)$$ 
for $\alpha$-stable L\'{e}vy motion.

The localisability of L\'{e}vy motion is a consequence of the fact that
it is $1/\alpha$-self-similar with stationary increments \cite{Fal6}. Its multistable
version and incremental moments are described in the following theorem:

\begin{theo} {\it (symmetric multistable L\'{e}vy motion)}.\label{esplevy}
Let $\alpha: [0,1] \to [c,d] \subset (1,2)$ be continuously differentiable. Let $(\Gamma_i)_{i \geq 1}$ be a sequence of arrival times of a Poisson process with unit arrival time,
$(V_i)_{i \geq 1}$ be a sequence of i.i.d. random variables with distribution $\hat{m}(dx) = dx$ on $[0,1]$, 
and $(\gamma_i)_{i \geq 1}$ be a sequence of i.i.d. random variables with distribution $P(\gamma_i=1)=P(\gamma_i=-1)=1/2$. 
Assume finally that the three sequences $(\Gamma_i)_{i \geq 1}$, $(V_i)_{i \geq 1}$, and $(\gamma_i)_{i \geq 1}$ are independent and 
define 
\begin{equation} 
 X(t,u)=C_{\alpha(u)}^{1/ \alpha(u)} \sum_{i=1}^{\infty} \gamma_i \Gamma_i^{-1/\alpha(u)}  \mathbf{1}_{[0,t]}(V_i)
 \end{equation}
and the symmetric multistable L\'{e}vy motion
\begin{displaymath}
Y(t) = X(t,t).
\end{displaymath}

Then for all $t \in (0,1)$ and $\eta < c$, when $\varep$ tends to 0,
 \begin{displaymath}
 \E \left[ |Y(t+\varep)-Y(t)|^{\eta}\right] \sim \frac{2^{\eta-1} \Gamma(1-\frac{\eta}{\alpha(t)})}{\eta \int_{0}^{\infty} u^{-\eta-1} \sin^2(u)du}\hspace{0.1cm}\varep^{\frac{\eta}{\alpha(t)}} .
\end{displaymath}
\end{theo}

\noindent
{\bf Proof} 

\noindent
See section {\ref{proo}}.

\begin{theo}\label{expolevy}
Let $Y$ be the symmetric multistable L\'{e}vy motion defined on $(0,1)$ with $\alpha: [0,1] \to [c,d] \subset (0,2)$. For all $t \in (0,1)$, almost surely,
\begin{displaymath}
\h \leq \frac{1}{\alpha(t)}.
\end{displaymath}
\end{theo}

\noindent
{\bf Proof} 

\noindent
See section {\ref{proo}}.

\begin{theo}\label{expolevyexact}
Let $u \in U \subset (0,1)$.
\begin{enumerate}
	\item If $0< \alpha(u) <1$, almost surely,
\begin{displaymath}
\mathcal{H}_u = \min\left(\frac{1}{\alpha(u)},\mathcal{H}_u^{\alpha}\right),
\end{displaymath}
where $\mathcal{H}_u^{\alpha}$ denotes the Hölder exponent of $\alpha$ at $u$,
at least when $\frac{1}{\alpha(u)} \neq \mathcal{H}_u^{\alpha}$.
	\item If $1 \leq \alpha(u) <2$, and  $\alpha$ is $\mathcal{C}^1$, almost surely,
\begin{displaymath}
\mathcal{H}_u = \frac{1}{\alpha(u)}.
\end{displaymath}
\end{enumerate}

\end{theo}

\noindent
{\bf Proof} 

\noindent
See section {\ref{proolevy}}.

\bigskip

Thus, in the case $0< \alpha(u) <1$, the regularity of the multistable Lévy motion is the smallest
number between $\frac{1}{\alpha(u)}$ and the regularity of the function $\alpha$ at $u$. This is very similar to
the case of the multifractional Brownian motion, where the Hölder exponent is the minimum between
the functional parameter $h$ and its regularity \cite{EH,EHJLV}. We conjecture that the same result 
holds even when $\alpha\geq1$. 

\section{Intermediate results}\label{IR}

Let $\varphi_X$ denote the characteristic function of the random variable $X$.
We first state the following almost obvious fact :

\begin{prop}\label{probinf}
Assume that for a given $t \in \bbbr$ there exists $\varep_0 > 0 $ such that 
\begin{displaymath}
\sup_{r \in B(0,\varep_0)} \int_{0}^{+\infty} \left| \varphi_{\frac{Y(t+r)-Y(t)}{r^{h(t)}}}(v)\right| dv <  + \infty ,
\end{displaymath}
where $Y$ is a symmetrical process. Then there exists $K >0$ which depends only on $t$ and $\varep_0$ such that for all $x >0$, and all $r \in (0, \varep_0)$,
\begin{displaymath}
\P \left( |Y(t+r)-Y(t)| < x \right) \leq K \frac{x}{r^{h(t)}}.
\end{displaymath}
If furthermore we suppose that $\sup_{t \in U} \sup_{r \in B(0,\varep_0)} \int_{0}^{+\infty} \left| \varphi_{\frac{Y(t+r)-Y(t)}{r^{h(t)}}}(v) \right| dv <  + \infty $, then for all $t \in U$, for all $r \in (0, \varep_0)$, $\P \left( |Y(t+r)-Y(t)| < x \right) \leq K \frac{x}{r^{h(t)}} $.
\end{prop}

\noindent
{\bf Proof} 

\noindent
This is a straightforward consequence of the inversion formula.  Let $x >0$ and $r < \varep_0$. Since $Y$ is a symmetrical process, $\varphi_{Y(t+r)-Y(t)}$ is an even function and
\begin{eqnarray*}
\P \left( |Y(t+r)-Y(t)| < x \right) & = &  \frac{1}{\pi} \left| \int_{0}^{+\infty} \varphi_{Y(t+r)-Y(t)} \left(\frac{v}{r^{h(t)}}\right) \sin \left( \frac{vx}{r^{h(t)}}\right) \frac{dv}{v} \right|\\
& \leq & \frac{1}{\pi} \frac{x}{r^{h(t)}} \sup_{r \in B(0,\varep_0)} \int_{0}^{+\infty} \left| \varphi_{\frac{Y(t+r)-Y(t)}{r^{h(t)}}}(v) \right| dv  \\
& \leq & K \frac{x}{r^{h(t)}}\Box \\
\end{eqnarray*}

We now consider multistable processes, first in the finite measure space case,
and then in the $\sigma$-finite measure space case:

\begin{prop}\label{msspfm}
Assuming (C1), (C2) and (C3), there exists $K_U >0$ such that for all $u \in U$, $v \in U$ and $x>0$,
\begin{displaymath}
\P \left(|X(v,v)-X(v,u)| > x \right) \leq K_U \left( \frac{|v-u|^d}{x^d}(1 +|\log \frac{|v-u |}{x}|^d) + \frac{|v-u|^c}{x^c}(1+|\log \frac{|v-u |}{x}|^c) \right) \\ .
\end{displaymath}
\end{prop}

\noindent
{\bf Proof} 

\noindent
See section \ref{proo}.

In the $\sigma$-finite space case, a similar property holds:

\begin{prop}\label{msspfm2}
Assuming (C1), (Cs2), (Cs3) and (Cs4), there exists $K_U >0$ such that for all $u \in U$, $v \in U$ and $x>0$,
\begin{displaymath}
\P \left(|X(v,v)-X(v,u)| > x \right) \leq K_U \left( \frac{|v-u|^d}{x^d}(1 +|\log \frac{|v-u |}{x}|^d) + \frac{|v-u|^c}{x^c}(1+|\log \frac{|v-u |}{x}|^c) \right) \\ .
\end{displaymath}
\end{prop}

\medskip

\noindent
{\bf Proof}

\noindent
We shall apply Proposition \ref{msspfm} to the function $g(t,w,x)=r(x)^{1 / \alpha(w)} f(t,w,x)$ on $(E, \mathcal{E}, \hat m )$.
\begin{itemize}
 \item By (C1), the family of functions $v \to f(t,v,x)$ is differentiable for all $(v,t)$ in $U^2$ and almost all $x$ in $E$ thus $v \to g(t,v,x)$  is differentiable too i.e (C1) holds for $g$.
\item Choose $\delta> \frac{d}{c} - 1$ such that (Cs2) holds. 
\begin{displaymath}
\sup_{w \in U} (|g(t,w,x)|^{\alpha(w)}) = r(x) \sup_{w \in U} (|f(t,w,x)|^{\alpha(w)})
\end{displaymath}
One has
\begin{eqnarray*}
\int_\bbbr \left[ \sup_{w \in U} (|g(t,w,x)|^{\alpha(w)}) \right]^{1+\delta} \hat m(dx) & =  & \int_\bbbr r(x)^{1+\delta} \left[\sup_{w \in U} (|f(t,w,x)|^{\alpha(w)})\right]^{1+\delta}  \hat m(dx)\\
 & = & \int_\bbbr \left[\sup_{w \in U} (|f(t,w,x)|^{\alpha(w)})\right]^{1+\delta} r(x)^{\delta}  m(dx)\\
\end{eqnarray*}
thus (C2) holds.

\item Choose $\delta> \frac{d}{c} - 1 $ such that (Cs3) and (Cs4) hold.
\begin{displaymath}
g'_u(t,w,x) = r(x)^{1/ \alpha(w)} (f'_u(t,w,x)-\frac{\alpha'(w)}{\alpha^2(w)}\log(r(x))f(t,w,x))
\end{displaymath}
and
\begin{displaymath}
\int_\bbbr \left[ \sup_{w \in U} (|g'_u(t,w,x)|^{\alpha(w)})\right]^{1+\delta} \hspace{0.1cm} \hat m(dx) 
\end{displaymath}
\begin{displaymath}
\leq \int_\bbbr \left[\sup_{w \in U} \left[ |f'_u(t,w,x)-\frac{\alpha'(w)}{\alpha^2(w)}\log(r(x))f(t,w,x)|^{\alpha(w)} \right] \right]^{1+\delta} r(x)^{\delta} m(dx).
\end{displaymath}
The inequality $|a+b|^\delta \leq \max(1, 2^{\delta-1})(|a|^\delta + |b|^\delta)$ shows that (C3) holds. 

Proposition \ref{msspfm} allows to conclude.
\end{itemize}

\Box

\begin{prop}\label{intfonc}
We suppose that there exists a function $h$ defined on $U$ such that (C8), (C10) and (C14) hold. Assuming (C1), (C6), (C7), (C11), (C12), (C13), (C15),  one has:

\begin{displaymath}
\sup_{r \in B(0,\varep)} \int_{0}^{+\infty} \varphi_{\frac{Y(t+r)-Y(t)}{r^{h(t)}}}(v) dv < +\infty.
\end{displaymath}

If in addition we suppose (Cu8), (Cu10), (Cu11), (Cu12),  (Cu14) and (Cu15), then
\begin{displaymath}
\sup_{t \in U} \sup_{r \in B(0,\varep)} \int_{0}^{+\infty} \varphi_{\frac{Y(t+r)-Y(t)}{r^{h(t)}}}(v) dv < +\infty.
\end{displaymath}
\end{prop}

\noindent
{\bf Proof} 

\noindent
The expression of the characteristic function $\varphi_{\frac{Y(t+r)-Y(t)}{r^{h(t)}}}$ is given in \cite{LGLV} :
\begin{displaymath}
\varphi_{\frac{Y(t+r)-Y(t)}{r^{h(t)}}}(v) = \exp \left( - 2 \int_{\bbbr} \int_{0}^{+\infty} \sin^2 \left( \frac{vC_{\alpha(t+r)}^{1/ \alpha(t+r)}f(t+r,t+r,x)}{2r^{h(t)}y^{1/ \alpha(t+r)}} - \frac{v C_{\alpha(t)}^{1/ \alpha(t)} f(t,t,x)}{2r^{h(t)}y^{1/ \alpha(t)}}\right)  dy \hspace*{0.1cm} m(dx)\right).
\end{displaymath}

For $v \leq 1$, $ \varphi_{\frac{Y(t+r)-Y(t)}{r^{h(t)}}}(v) \leq 1$. For $v \geq 1$, we fix $\varep < \frac{1}{d}$. Lemma (\ref{lemsin}) entails that there exists $K_U >0$ such that 
\begin{displaymath}
\varphi_{\frac{Y(t+r)-Y(t)}{r^{h(t)}}}(v) \leq \exp \left( - \int_{\bbbr} \int_{\frac{K_U v^{\frac{d}{1-\varep d}}}{r}} \left| \frac{v C_{\alpha(t+r)}^{1/ \alpha(t+r)}f(t+r,t+r,x)}{2r^{h(t)}y^{1/ \alpha(t+r)}}  - \frac{v C_{\alpha(t)}^{1/ \alpha(t)}f(t,t,x)}{2r^{h(t)}y^{1/ \alpha(t)}} \right|^2 dy \hspace*{0.1cm} m(dx) \right).
\end{displaymath}

Let 
\begin{displaymath}
 N(v,t,r) = \int_{\bbbr} \int_{\frac{K_U v^{\frac{d}{1-\varep d}}}{r}} \left| \frac{v C_{\alpha(t+r)}^{1/ \alpha(t+r)}f(t+r,t+r,x)}{2r^{h(t)}y^{1/ \alpha(t+r)}}  - \frac{v C_{\alpha(t)}^{1/ \alpha(t)}f(t,t,x)}{2r^{h(t)}y^{1/ \alpha(t)}} \right|^2 dy \hspace*{0.1cm} m(dx).
\end{displaymath}
Using lemma ($\ref{borneN}$), there exist $K_U >0$ and $\varep_0 > 0$ such that for all $v \geq 1$, 
\begin{displaymath}
 N(v,t,r) \geq K_U v^{2+\frac{d}{1-\varep d}(1- \frac{2}{c})}.
\end{displaymath}
The inequality becomes
\begin{displaymath}
\varphi_{\frac{Y(t+r)-Y(t)}{r^{h(t)}}}(v) \leq \exp \left( - K_U v^{2+\frac{d}{1-\varep d}(1- \frac{2}{c})} \right),
\end{displaymath}
and
\begin{eqnarray*}
\int_{0}^{+\infty} \varphi_{\frac{Y(t+r)-Y(t)}{r^{h(t)}}}(v) dv & \leq & 1 + \int_{1}^{\infty} \exp \left( - K_U v^{2+\frac{d}{1-\varep d}(1- \frac{2}{c})} \right) dv \\
& < & + \infty\\
\end{eqnarray*}
\Box

\section{Proofs and technical results}\label{proo}

{\bf Proof of proposition \ref{msspfm}}

We proceed as in \cite{LGLV}. Note that condition {\bf (C2)} implies that there exists $\delta > \frac{d}{c}-1$ such that :
\begin{equation}\label{kercond2f}
\sup_{t \in U}  \int_\bbbr \left[ \sup_{w \in U} \left[\left|f(t,w,x)\log|f(t,w,x)|\right|^{\alpha(w)}\right] \right]^{1+\delta} \hspace{0.1cm} \hat m(dx) < \infty.
\end{equation}
The function $u \mapsto C^{1/\alpha(u)}_{\alpha(u)}$ is a $C^1$ function since $\alpha(u)$ ranges in $[c,d] \subset (0,2)$. We shall denote $a(u) = (m(E))^{1/\alpha(u)} C^{1/\alpha(u)}_{\alpha(u)}$. The function $a$ is thus also $C^1$. Let $(u,v) \in U^2$.We estimate:

$$X(v,v)-X(v,u) = \sum_{i=1}^{\infty} \gamma_i (\Phi_i(v)-\Phi_i(u)) + \sum_{i=1}^{\infty}\gamma_i (\Psi_i(v)-\Psi_i(u)),$$

where
 \begin{displaymath}
\Phi_i(u) = a(u) i^{-1/\alpha(u)} f(v,u,V_i)
\end{displaymath}
and 
\begin{displaymath}
\Psi_i(u)=a(u) \left( \Gamma_i^{-1/\alpha(u)} - i^{-1/\alpha(u)} \right) f(v,u,V_i). 
\end{displaymath}
Thanks to the assumptions on $a $ and $f$, $\Phi_i$ and $\Psi_i$ are differentiable and one computes :
$$\Phi_i'(u) = a'(u) i^{-1/\alpha(u)} f(v,u,V_i) + a(u)i^{-1/\alpha(u)} f'_u(v,u,V_i) + a(u)\frac{\alpha'(u)}{\alpha(u)^2} \log(i) i^{-1/\alpha(u)} f(v,u,V_i),$$
and
$$
\Psi_i'(u) = a'(u) \left( \Gamma_i^{-1/\alpha(u)} - i^{-1/\alpha(u)} \right) f(v,u,V_i) + a(u)\left( \Gamma_i^{-1/\alpha(u)} - i^{-1/\alpha(u)} \right) f'_u(v,u,V_i)
$$
$$+ a(u)\frac{\alpha'(u)}{\alpha(u)^2} \left( \log(\Gamma_i) \Gamma_i^{-1/\alpha(u)} - \log(i) i^{-1/\alpha(u)} \right) f(v,u,V_i).
$$

Using the mean value theorem, there exists a sequence of independent random numbers $w_i \in (u,v)$ (or $(v,u)$) and a sequence of random numbers $x_i \in (u,v)$ (or $(v,u)$) such that:
\begin{equation}\label{decompchamp}
X(v,u)-X(v,v) = (u-v)\sum_{i=1}^{\infty} (Z_i^1+Z_i^2+Z_i^3) + (u-v)\sum_{i=1}^{\infty} (Y_i^1+Y_i^2+Y_i^3),
\end{equation}
where
$$Z_i^1 = \gamma_i a'(w_i)i^{-1/\alpha(w_i)} f(v,w_i,V_i),$$
$$Z_i^2 = \gamma_i a(w_i)i^{-1/\alpha(w_i)} f'_u(v,w_i,V_i),$$
$$Z_i^3 = \gamma_i a(w_i)\frac{\alpha'(w_i)}{\alpha(w_i)^2} \log(i) i^{-1/\alpha(w_i)} f(v,w_i,V_i),$$
$$Y_i^1 = \gamma_i a'(x_i)\left( \Gamma_i^{-1/\alpha(x_i)} - i^{-1/\alpha(x_i)} \right) f(v,x_i,V_i),$$
$$Y_i^2 = \gamma_i a(x_i)\left( \Gamma_i^{-1/\alpha(x_i)} - i^{-1/\alpha(x_i)} \right) f'_u(v,x_i,V_i),$$
$$Y_i^3 = \gamma_i a(x_i)\frac{\alpha'(x_i)}{\alpha(x_i)^2} \left( \log(\Gamma_i) \Gamma_i^{-1/\alpha(x_i)} - \log(i) i^{-1/\alpha(x_i)} \right) f(v,x_i,V_i).$$
Note that each $w_i$ depends on $a,f,\alpha,u,v,V_i$, and each $x_i$ depends on $a,f,\alpha,u,v,V_i, \Gamma_i$ but not on $\gamma_i$. This remark will be useful in the sequel.

In \cite{LGLV}, it is proved that each series $\sum_{i=1}^{\infty}\limits Z_i^j$ and $\sum_{i=1}^{\infty}\limits Y_i^j$ , $j=1,2,3, $ converges almost surely. Let $x>0$. We consider $\P \left( | \sum_{i=1}^{\infty}\limits Z_i^j| > x \right)$ and $\P \left( |\sum_{i=1}^{\infty}\limits Y_i^j | > x \right)$ for $j=1,2,3$.
\smallskip

Let $\eta \in (0,\min(\frac{2c}{d}-1,\frac{c}{d}(\delta+1)-1))$.
Markov inequality yields
 \begin{eqnarray*}
 \P\left( | \sum_{i=1}^{\infty} Z_i^j| > x \right) & \leq & \frac{1}{x^d} \E \left[ |\sum_{i=1}^{\infty}  Z_i^j|^d \right]\\
 & \leq & \frac{1}{x^d} \left(\E \left[ |\sum_{i=1}^{\infty}  Z_i^j|^{d(1+\eta)} \right]\right)^{\frac{1}{1+\eta}}.\\
 \end{eqnarray*}
 The random variables $Z_i^j$ are independent with mean $0$ thus, by theorem 2 of \cite{BE}: 
 $$\E \left[ | \sum_{i=1}^{+\infty}\limits Z_i^j |^{d(1+\eta)} \right] \leq 2 \sum_{i=1}^{+\infty}\limits \E [ | Z_i^j|^{d(1+\eta)}].$$\par
 
 For $j=1$,
 \begin{eqnarray*}
 \E \left[|Z_i^1|^{d(1+\eta)}\right] & = & \E \left[ |a'(w_i)|^{d(1+\eta)} i^{-\frac{d(1+\eta)}{\alpha(w_i)}} |f(v,w_i,V_i)|^{d(1+\eta)} \right]\\
 & \leq & \frac{K_U}{i^{1+\eta}} \E \left[ \left( \sup_{w \in B(u,\varep)} |f(v,w,V_i) |^{\alpha(w)}\right)^{\frac{d(1+\eta)}{\alpha(w_i)}}\right]\\
 & \leq & \frac{K_U}{i^{1+\eta}} \E \left[ (\sup_{w \in B(u,\varep)} |f(v,w,V_1) |^{\alpha(w)})^{1+\eta} + (\sup_{w \in B(u,\varep)} |f(v,w,V_1) |^{\alpha(w)})^{\frac{d}{c}(1+\eta)}\right]\\
 & \leq & \frac{K_U}{i^{1+\eta}}.\\
 \end{eqnarray*}

 For $j=2$,
 \begin{eqnarray*}
 \E \left[|Z_i^2|^{d(1+\eta)}\right] & \leq & \frac{K_U}{i^{1+\eta}} \E \left[ (\sup_{w \in B(u,\varep)} |f'_u(v,w,V_1) |^{\alpha(w)})^{1+\eta} + (\sup_{w \in B(u,\varep)} |f'_u(v,w,V_1) |^{\alpha(w)})^{\frac{d}{c}(1+\eta)}\right]\\
 & \leq & \frac{K_U}{i^{1+\eta}}.\\
 \end{eqnarray*}
 For $j=3$,
 \begin{eqnarray*}
   \E \left[|Z_i^3|^{d(1+\eta)}\right] & = & \E\left[ |a(w_i)\frac{\alpha'(w_i)}{\alpha(w_i)^2} |^{d(1+\eta)} |f(v,w_i,V_i) |^{d(1+\eta)} \frac{(\log i)^{d(1+\eta)}}{i^{\frac{d(1+\eta)}{\alpha(w_i)}}}\right] \\
   & \leq & K_U \frac{(\log i )^{d(1+\eta)}}{i^{1+\eta}}.
 \end{eqnarray*}
 Finally, $\sup_{v \in U}\limits \sum_{i=1}^{+\infty}\limits \E \left[ |Z_i^j|^{d(1+\eta)}\right] < +\infty$ thus 
 
 \begin{displaymath}
 \P\left( | \sum_{i=1}^{\infty} Z_i^j| > x \right) \leq \frac{K_U}{x^d}.
 \end{displaymath}
We consider now $\P \left( |\sum_{i=1}^{\infty}\limits Y_i^j | > x \right)$ for $j=1,2,3$.

\begin{displaymath}
\P\left( | \sum_{i=1}^{\infty} Y_i^j| > x \right)  \leq  \P\left( |Y_1^j| \geq \frac{x}{2} \right)+\P\left( | \sum_{i=2}^{\infty} Y_i^j| \geq \frac{x}{2} \right).
\end{displaymath}

Since $\P\left( | \sum_{i=2}^{\infty}\limits Y_i^j| \geq \frac{x}{2} \right) \leq \frac{2^d}{x^d} \left( \E \left[ | \sum_{i=2}^{\infty}\limits Y_i^j|^{d(1+\eta)} \right]\right)^{\frac{1}{1+\eta}}$, we want to apply theorem 2 of \cite{BE} again. 
Let $S_m=\sum_{i=1}^{m}\limits Y_i^j$ and write $Y_i^j=\gamma_i W_i^j$. Note that $\gamma_i$ is independent of $W_i^j$ and $S_{i-1}$.
\begin{eqnarray*}
\E \left( Y_{m+1}^j | S_m \right) & = & \E \left( \E (Y_{m+1}^j | S_m,W_{m+1}) | S_m \right) \\
& = & \E \left( \E (\gamma_{m+1} W_{m+1}^j | S_m,W_{m+1}) | S_m \right) \\
& = & \E \left( W_{m+1}^j \E (\gamma_{m+1} | S_m,W_{m+1}) | S_m \right) \\
& = & \E \left( W_{m+1}^j \E (\gamma_{m+1}) | S_m \right) \\
& = & 0.\\
\end{eqnarray*}
We apply theorem 2 of \cite{BE} with $(d(1+\eta)<2)$,
\begin{displaymath}
\E \left[ | \sum_{i=2}^{\infty}Y_i^j|^{d(1+\eta)} \right] \leq 2 \sum_{i=2}^{\infty} \E |Y_i^j|^{d(1+\eta)},
\end{displaymath}
and
\begin{displaymath}
\P\left( | \sum_{i=1}^{\infty} Y_i^j| > x \right)  \leq  \P\left( |Y_1^j| \geq \frac{x}{2} \right) + \frac{2^{d}}{x^d} \left( 2 \sum_{i=2}^{\infty} \E |Y_i^j|^{d(1+\eta)}\right)^{\frac{1}{1+\eta}}.
\end{displaymath}
For $j=1$,

\begin{eqnarray*}
\P\left( |Y_1^1| \geq \frac{x}{2} \right) & = & \P \left( |a'(x_1)|^{\alpha(x_1)} | \frac{1}{\Gamma_1^{1/ \alpha(x_1)}}-1 |^{\alpha(x_1)} |f(v,x_1,V_1) |^{\alpha(x_1)} \geq \frac{x^{\alpha(x_1)}}{2^{\alpha(x_1)}}\right)\\
& \leq & \P \left( | \frac{1}{\Gamma_1^{1/ \alpha(x_1)}}-1 |^{\alpha(x_1)} \sup_{w \in B(u,\varep)} |f(v,w,V_1)|^{\alpha(w)} \geq K_U x^{\alpha(x_1)}\right).\\
\end{eqnarray*}

For $x <1$,

\begin{eqnarray*}
\P\left( |Y_1^1| \geq \frac{x}{2} \right) & \leq & \P \left( | \frac{1}{\Gamma_1^{1/ \alpha(x_1)}}-1 |^{\alpha(x_1)} \sup_{w \in B(u,\varep)} |f(v,w,V_1)|^{\alpha(w)} \geq K_U x^{d}\right)\\
& \leq & \P \left( \{ \sup_{w \in B(u,\varep)} |f(v,w,V_1)|^{\alpha(w)} \geq K_U x^{d}\} \cap \{\Gamma_1 > 1 \}\right)\\
& + & \P \left( \{ | \frac{1}{\Gamma_1} |\Gamma_1^{1/ \alpha(x_1)}-1 |^{\alpha(x_1)} \sup_{w \in B(u,\varep)} |f(v,w,V_1)|^{\alpha(w)} \geq K_U x^{d} \} \cap \{ \Gamma_1 \leq 1\}\right).\\
\end{eqnarray*}

\begin{eqnarray*}
 \P \left( \{ \sup_{w \in B(u,\varep)} |f(v,w,V_1)|^{\alpha(w)} \geq K_U x^{d}\} \cap \{\Gamma_1 > 1 \}\right) & \leq & \frac{K_U}{x^d} \E \left( \sup_{w \in B(u,\varep)} |f(v,w,V_1)|^{\alpha(w)}\right)\\
 & \leq & \frac{K_U}{x^d}.\\ 
\end{eqnarray*}
Let $W(v,x)= \sup_{w \in B(u,\varep)} |f(v,w,x)|^{\alpha(w)}$ and $F_{v,V_1}$ be the distribution of $W(v,V_1)$.
\begin{eqnarray*}
 \P \left( \{ | \frac{1}{\Gamma_1} |\Gamma_1^{1/ \alpha(x_1)}-1 |^{\alpha(x_1)} W(v,V_1) \geq K_U x^{d} \} \cap \{ \Gamma_1 \leq 1\}\right) & \leq  & \P \left(   W(v,V_1) \geq K_U x^{d} \Gamma_1 \right)\\
 & = & \int_{0}^{+\infty} \P \left( z \geq K_U x^{d} \Gamma_1 \right) F_{v,V_1}(dz) \\
 & = & \int_{0}^{+\infty} \left( 1 - e^{-\frac{z}{K_U x^d}}\right) F_{v,V_1}(dz) \\
 & \leq & \int_{0}^{+\infty} \frac{z}{K_U x^d}F_{v,V_1}(dz) \\
 & \leq & \frac{K_U}{ x^d}.
\end{eqnarray*}

For $x \geq 1$,
\begin{eqnarray*}
 \P\left( |Y_1^1| \geq \frac{x}{2} \right) & \leq & \P \left( | \frac{1}{\Gamma_1^{1/ \alpha(x_1)}}-1 |^{\alpha(x_1)} \sup_{w \in B(u,\varep)} |f(v,w,V_1)|^{\alpha(w)} \geq K_U x^{c}\right)\\
 \P\left( |Y_1^1| \geq \frac{x}{2} \right) & \leq & \frac{K_U}{x^c}.\\
\end{eqnarray*}

For $i \geq 2$,
\begin{eqnarray*}
 \E |Y_i^1|^{d(1+\eta)} & = & \E \left( |a'(x_i)|^{d(1+\eta)} | \Gamma_i^{-1/\alpha(x_i)} - i^{-1/\alpha(x_i)}|^{d(1+\eta)} \left(|f(v,x_i,V_i)|^{\alpha(x_i)} \right)^{\frac{d(1+\eta)}{\alpha(x_i)}}\right)\\
 & \leq & K_U \E \left( i^{-\frac{d(1+\eta)}{\alpha(x_i)}} W(v,V_i)^{\frac{d(1+\eta)}{\alpha(x_i)}} |(\frac{i}{\Gamma_i} )^{1/ \alpha(x_i)} - 1 |^{d(1+\eta)} \right)\\
 & \leq & \frac{K_U}{i^{1+\eta}} \E \left( \left[  W(v,V_i)^{1+\eta} +  W(v,V_i)^{\frac{d}{c}(1+\eta) }\right]\left[ |(\frac{i}{\Gamma_i} )^{1/ c} - 1 |^{d(1+\eta)} + |(\frac{i}{\Gamma_i} )^{1/d} - 1 |^{d(1+\eta)}\right]\right)\\
  & \leq & \frac{K_U}{i^{1+\eta}} \E \left(  W(v,V_i)^{1+\eta} +  W(v,V_i)^{\frac{d}{c}(1+\eta) } \right) \E \left( |(\frac{i}{\Gamma_i} )^{1/ c} - 1 |^{d(1+\eta)} + |(\frac{i}{\Gamma_i} )^{1/d} - 1 |^{d(1+\eta)}\right).\\
\end{eqnarray*}
Using the fact that $\eta \leq \delta$ and $\frac{d}{c}(1+\eta) \leq 1+\delta$,

\begin{eqnarray*}
\E \left(  W(v,V_i)^{1+\eta} +  W(v,V_i)^{\frac{d}{c}(1+\eta) } \right) & = & \E \left(  W(v,V_1)^{1+\eta} +  W(v,V_1)^{\frac{d}{c}(1+\eta) } \right)\\
& \leq & K_U,
\end{eqnarray*}
\begin{eqnarray*}
\E |(\frac{i}{\Gamma_i} )^{1/ c} - 1 |^{d(1+\eta)} & \leq  & K_U(1+ \E \left( (\frac{i}{\Gamma_i})^{\frac{d}{c}(1+\eta)}\right))\\
& \leq & K_U,
\end{eqnarray*}
and
\begin{displaymath}
\E |(\frac{i}{\Gamma_i} )^{1/ d} - 1 |^{d(1+\eta)} \leq K_U.
\end{displaymath}
As a consequence:

\begin{displaymath}
\sup_{v \in U}\limits \sum_{i=2}^{+\infty} \E |Y_i^1|^{d(1+\eta)}  \leq K_U
\end{displaymath}
and
\begin{displaymath}
\P \left( \left| \sum_{i=1}^{\infty} Y_i^1 \right| > x \right) \leq K_U \left( \frac{1}{x^c} + \frac{1}{x^d} \right).
\end{displaymath}

For $j=2$, since the conditions required on $(a',f)$ are also satisfied by $(a,f'u)$, one
gets in a similar fashion
\begin{displaymath}
\P \left( \left| \sum_{i=1}^{\infty} Y_i^2 \right| > x \right) \leq K_U \left( \frac{1}{x^c} + \frac{1}{x^d} \right).
\end{displaymath}

For $j=3$,

\begin{eqnarray*}
\P \left( |Y_1^3 | \geq \frac{x}{2}  \right) & = & \P \left( |a(x_1) \frac{\alpha'(x_1)}{\alpha(x_1)^2}\log(\Gamma_1) \Gamma_1^{-1/\alpha(x_1)} f(v,x_1,V_1) | \geq \frac{x}{2} \right)\\
& \leq & \P \left( K_U \frac{|f(v,x_1,V_1)|^{\alpha(x_1)}}{x^{\alpha(x_1)}} \geq \frac{\Gamma_1}{| \log \Gamma_1|^{\alpha(x_1)}} \right). \\
\end{eqnarray*}

Let $g(z)= \frac{z}{| \log z |^{\alpha(x_1)}}$, for $z < 1$.

$g$ is a one-to-one increasing function, and for all $z < 1$ such that $z|\log z|^{\alpha(x_1)} <1$ and $|1+\alpha(x_1) \frac{\log |\log z|}{|\log z|}|^{\alpha(x_1)} \leq 2$,

 \begin{displaymath}
 g\left( z|\log z|^{\alpha(x_1)}\right) = \frac{z|\log z|^{\alpha(x_1)}}{|\log z +  \alpha(x_1) \log | \log z| |^{\alpha(x_1)}} \geq \frac{z}{2}
 \end{displaymath}
thus $g^{-1}(\frac{z}{2}) \leq z |\log z|^{\alpha(x_1)}$.
\newline
Fix $A>0$ such that for all $0<z<A$, $g^{-1}(z) \leq 2z |\log 2 + \log z|^{\alpha(x_1)} $ {\it i.e.} 
\begin{displaymath}
g^{-1}(z) \leq K_U z |\log z|^{\alpha(x_1)}.
\end{displaymath}

Let $B=\{  K_U \frac{|f(v,x_1,V_1)|^{\alpha(x_1)}}{x^{\alpha(x_1)}} \geq \frac{\Gamma_1}{| \log \Gamma_1|^{\alpha(x_1)}} \}$.
\begin{displaymath}
\P (B)= \P(B \cap \{\Gamma_1 >1 \} ) +\P(B \cap \{\Gamma_1 <1 \} \cap \{ 0 \leq K_U \frac{|f(v,x_1,V_1)|^{\alpha(x_1)}}{x^{\alpha(x_1)}} \leq A\} ) 
\end{displaymath}
\begin{displaymath}
+ \P (B \cap \{ \Gamma_1 <1\} \cap \{K_U \frac{|f(v,x_1,V_1)|^{\alpha(x_1)}}{x^{\alpha(x_1)}} > A \}).
\end{displaymath}

Each of these three terms will be treated separately.

\begin{eqnarray*}
\bullet \hspace*{0.3cm}\P \left(B \cap \{\Gamma_1 >1 \} \right) & \leq & \P \left(K_U \frac{|f(v,x_1,V_1)|^{\alpha(x_1)}}{x^{\alpha(x_1)}} |\log \Gamma_1 |^{\alpha(x_1)} \geq 1\right)\\
& \leq & \P \left( K_U \sup_{w \in B(u,\varep)} |f(v,w,V_1)|^{\alpha(w)} ( |\log \Gamma_1 |^{c} + |\log \Gamma_1 |^{d}) \geq x^{\alpha(x_1)} \right) .\\
\end{eqnarray*}
For $x \geq 1$,
\begin{eqnarray*}
\P \left(B \cap \{\Gamma_1 >1 \} \right) & \leq & \P \left( K_U \sup_{w \in B(u,\varep)} |f(v,w,V_1)|^{\alpha(w)} ( |\log \Gamma_1 |^{c} + |\log \Gamma_1 |^{d}) \geq x^{c} \right)\\
& \leq & \frac{K_U}{x^c} \E \left( \sup_{w \in B(u,\varep)} |f(v,w,V_1)|^{\alpha(w)}\right) \E \left(|\log \Gamma_1 |^{c} + |\log \Gamma_1 |^{d} \right)\\
& \leq & \frac{K_U}{x^c}.
\end{eqnarray*}
For $x < 1$,
\begin{eqnarray*}
\P\left(B \cap \{\Gamma_1 >1 \} \right) & \leq & \P \left( K_U \sup_{w \in B(u,\varep)} |f(v,w,V_1)|^{\alpha(w)} ( |\log \Gamma_1 |^{c} + |\log \Gamma_1 |^{d}) \geq x^{d} \right)\\
& \leq & \frac{K_U}{x^d}.
\end{eqnarray*}
\begin{eqnarray*}
\bullet \hspace*{0.3cm}\P \left(B \cap \{ \Gamma_1 <1\} \cap \{K_U \frac{|f(v,x_1,V_1)|^{\alpha(x_1)}}{x^{\alpha(x_1)}} > A \}\right) & \leq & \P \left(K_U |f(v,x_1,V_1)|^{\alpha(x_1)} \geq A x^{\alpha(x_1)}\right)\\
& \leq & \frac{K_U}{x^c} +\frac{K_U}{x^d}.\\
\end{eqnarray*}

\begin{eqnarray*}
\bullet \hspace*{0.3cm}\P \left(B \cap \{ \Gamma_1 <1\} \cap \{0 \leq K_U \frac{|f(v,x_1,V_1)|^{\alpha(x_1)}}{x^{\alpha(x_1)}} \leq  A \}\right) &  & \\
\end{eqnarray*}
\begin{eqnarray*}
& = & \P \left(\{ g(\Gamma_1) \leq K_U \frac{|f(v,x_1,V_1)|^{\alpha(x_1)}}{x^{\alpha(x_1)}}\}  \cap \{ \Gamma_1 <1\} \cap \{0 \leq K_U \frac{|f(v,x_1,V_1)|^{\alpha(x_1)}}{x^{\alpha(x_1)}} \leq  A \}\right)\\
& \leq & \P\left( \Gamma_1 \leq K_U \frac{|f(v,x_1,V_1)|^{\alpha(x_1)}}{x^{\alpha(x_1)}} + K_U \frac{|f(v,x_1,V_1)|^{\alpha(x_1)}}{x^{\alpha(x_1)}}\left| \log  \frac{|f(v,x_1,V_1)|^{\alpha(x_1)}}{x^{\alpha(x_1)}} \right|^{\alpha(x_1)}\right)\\
& \leq & \P\left(  \Gamma_1 \leq  K_U |f(v,x_1,V_1)|^{\alpha(x_1)} ( \frac{1 +| \log x|^{\alpha(x_1)}}{x^{\alpha(x_1)}}) +  K_U \frac{| |f(v,x_1,V_1)| \log |f(v,x_1,V_1)||^{\alpha(x_1)}}{x^{\alpha(x_1)}}\right).\\
\end{eqnarray*}
With $ W(v,x)= \sup_{w \in B(u,\varep)} |f(v,w,x)|^{\alpha(w)}$ and $Z(v,x)= \sup_{w \in B(u,\varep)} |f(v,w,x)\log | f(v,w,x)| |^{\alpha(w)} $,
\begin{eqnarray*}
\P \left(B \cap \{ \Gamma_1 <1\} \cap \{0 \leq K_U \frac{|f(v,x_1,V_1)|^{\alpha(x_1)}}{x^{\alpha(x_1)}} \leq  A \}\right) &  &\\
\end{eqnarray*}
\begin{eqnarray*}
& \leq &\P\left(  \Gamma_1 \leq  K_U   W(v,V_1) (\frac{1 +| \log x|^{\alpha(x_1)}}{x^{\alpha(x_1)}}) +  K_U \frac{Z(v,V_1)}{x^{\alpha(x_1)}}\right)\\
& \leq & \P\left(  \Gamma_1 \leq  K_U   W(v,V_1) (\frac{1 +| \log x|^{\alpha(x_1)}}{x^{\alpha(x_1)}}) \right) +\P\left(  \Gamma_1 \leq  K_U   Z(v,V_1) (\frac{1 +| \log x|^{\alpha(x_1)}}{x^{\alpha(x_1)}}) \right).\\
\end{eqnarray*}
Since $\frac{1 +| \log x|^{\alpha(x_1)}}{x^{\alpha(x_1)}} \leq  K_U(\frac{1 +| \log x|^{c}}{x^{c}}+\frac{1 +| \log x|^{d}}{x^{d}}) $,
\begin{eqnarray*}
\P\left(  \Gamma_1 \leq  K_U   W(v,V_1) (\frac{1 +| \log x|^{\alpha(x_1)}}{x^{\alpha(x_1)}}) \right) & \leq &\P\left(  \Gamma_1 \leq  K_U   W(v,V_1) (\frac{1 +| \log x|^{c}}{x^{c}}+\frac{1 +| \log x|^{d}}{x^{d}})\right)\\
& \leq & K_U(\frac{1 +| \log x|^{c}}{x^{c}}+\frac{1 +| \log x|^{d}}{x^{d}}),\\
\end{eqnarray*}
and
\begin{eqnarray*}
\P\left(  \Gamma_1 \leq  K_U   Z(v,V_1) (\frac{1 +| \log x|^{\alpha(x_1)}}{x^{\alpha(x_1)}}) \right) & \leq & \P\left(  \Gamma_1 \leq  K_U   Z(v,V_1) (\frac{1 +| \log x|^{c}}{x^{c}}+\frac{1 +| \log x|^{d}}{x^{d}})\right).\\
\end{eqnarray*}

Denoting $G_{v,V_1}$ the distribution of $Z(v,V_1)$,
\begin{eqnarray*}
\P\left(  \Gamma_1 \leq  K_U   Z(v,V_1) (\frac{1 +| \log x|^{c}}{x^{c}}+\frac{1 +| \log x|^{d}}{x^{d}})\right)\\
\end{eqnarray*}
\begin{eqnarray*}
 &= & \int_{0}^{+ \infty} (1 - \exp (-K_U (\frac{1 +| \log x|^{c}}{x^{c}}+\frac{1 +| \log x|^{d}}{x^{d}}) z )  ) G_{v,V_1}(dz) \\
& \leq & K_U (\frac{1 +| \log x|^{c}}{x^{c}}+\frac{1 +| \log x|^{d}}{x^{d}}) \int_{0}^{+ \infty} z G_{v,V_1}(dz)\\
& \leq & K_U (\frac{1 +| \log x|^{c}}{x^{c}}+\frac{1 +| \log x|^{d}}{x^{d}}),\\
\end{eqnarray*}
since $\sup_{v \in B(u,\varep)} \E \left( Z(v,V_1) \right)< +\infty$.
\newline

Finally,
\begin{displaymath}
\P \left(B \cap \{ \Gamma_1 <1\} \cap \{0 \leq K_U \frac{|f(v,x_1,V_1)|^{\alpha(x_1)}}{x^{\alpha(x_1)}} \leq  A \}\right) \leq  K_U (\frac{1 +| \log x|^{c}}{x^{c}}+\frac{1 +| \log x|^{d}}{x^{d}})
\end{displaymath}
and
\begin{displaymath}
 \P \left( |Y_1^3 | \geq \frac{x}{2}  \right) \leq  K_U (\frac{1 +| \log x|^{c}}{x^{c}}+\frac{1 +| \log x|^{d}}{x^{d}}).
\end{displaymath}

For $i \geq 2$,

\begin{eqnarray*}
\E |Y_i^3|^{d(1+\eta)} & \leq & K_U \frac{| \log i|^{d(1+\eta)}}{i^{1+\eta}} \E \left( W(v,V_i)^{1+\eta} + W(v,V_i)^{\frac{d}{c}(1+\eta)} \right) \E \left( |\frac{\log \Gamma_i}{\log i}||(\frac{i}{\Gamma_i})^{1/\alpha(x_i)} -1|^{d(1+\eta)}  \right) \\
 & \leq &  K_U \frac{| \log i|^{d(1+\eta)}}{i^{1+\eta}} \E \left( |\frac{\log \Gamma_i}{\log i}||(\frac{i}{\Gamma_i})^{1/\alpha(x_i)} -1|^{d(1+\eta)}  \right) \\
 & \leq & K_U \frac{| \log i|^{d(1+\eta)}}{i^{1+\eta}} \E \left( |\frac{\log \Gamma_i}{\log i}||(\frac{i}{\Gamma_i})^{1/c} -1|^{d(1+\eta)} + |\frac{\log \Gamma_i}{\log i}||(\frac{i}{\Gamma_i})^{1/d} -1|^{d(1+\eta)}\right)\\
 & \leq & K_U \frac{| \log i|^{d(1+\eta)}}{i^{1+\eta}},\\
\end{eqnarray*}
thus 
\begin{displaymath}
\sup_{v \in U}\limits \sum_{i=2}^{+\infty} \E |Y_i^3|^{d(1+\eta)}  \leq K_U
\end{displaymath}
and
\begin{displaymath}
\P \left( \left| \sum_{i=1}^{\infty} Y_i^3 \right| > x \right) \leq K_U (\frac{1 +| \log x|^{c}}{x^{c}}+\frac{1 +| \log x|^{d}}{x^{d}}).
\end{displaymath}

Let us go back to $\P \left(|X(v,v)-X(v,u)| > x \right) $.
\begin{eqnarray*}
 \P \left(|X(v,v)-X(v,u)| > x \right) & = & \P \left(|u-v||\sum_{i=1}^{\infty} (Z_i^1+Z_i^2+Z_i^3+Y_i^1+Y_i^2+Y_i^3) | > x \right) \\
 & \leq & \sum_{j=1}^{3} \left( \P\left( | \sum_{i=1}^{\infty} Z_i^j| \geq \frac{x}{6|u-v|} \right) + \P\left( | \sum_{i=1}^{\infty} Y_i^j| \geq \frac{x}{6|u-v|} \right) \right) \\
 & \leq & K_U \left( \frac{|v-u|^d}{x^d}(1+|\log \frac{|v-u |}{x}|^d) + \frac{|v-u|^c}{x^c}(1+|\log \frac{|v-u |}{x}|^c) \right) \\
\end{eqnarray*}

and the proof is complete \Box

\begin{lem} \label{lemN}
Assume (C11), (C12), (C14), (C15). There exists a function $l\geq 0$ such that
\begin{displaymath}
\lim_{r \rightarrow 0} |\Delta(r,t) - l(t)| =0,
\end{displaymath}
where 
\begin{displaymath}
\Delta(r,t) =: \frac{1}{r^{2h(t)}}\int_{\bbbr} \int_{\frac{K}{r}} \left| \frac{C_{\alpha(t+r)}^{1/ \alpha(t+r)}}{y^{1/ \alpha(t+r)}} f(t+r,t+r,x) - \frac{C_{\alpha(t)}^{1/ \alpha(t)}}{y^{1/ \alpha(t)}} f(t,t,x) \right|^2 dy \hspace*{0.1cm} m(dx).
\end{displaymath}
Assuming in addition (Cu11), (Cu12), (Cu14), (Cu15), the convergence is uniform on $U$.
\end{lem}

\noindent
{\bf Proof} 

\noindent
Let $l(t) = \frac{C_{\alpha(t)}^{2/ \alpha(t)} K^{1- \frac{2}{\alpha(t)}}}{\frac{2}{\alpha(t)}-1} g(t).$ Note that condition (C14) implies the following:
\begin{equation}\label{C16}
 \forall \varep > 0, \exists K_U >0 , \forall r \leq \varep, \frac{1}{|r|^{1+2(h(t)-\frac{1}{\alpha(t)})}} \int_{\bbbr} \left| f(t+r,t,x) - f(t,t,x) \right|^2 m(dx)\leq K_U.
\end{equation}
The uniform condition (Cu14) implies also that:
\begin{equation}\label{Cu16}
 \exists K_U >0 , \forall v \in U, \forall u \in U,  \frac{1}{|v-u|^{1+2(h(u)-\frac{1}{\alpha(u)})}} \int_{\bbbr} \left| f(v,u,x) - f(u,u,x) \right|^2 m(dx)\leq K_U.
\end{equation}

Expanding the square, we can write $\Delta(r,t) - l(t) = \Delta_1(r,t) + \Delta_2(r,t) + \Delta_3(r,t)$ where
\begin{displaymath}
\Delta_1(r,t) = \frac{1}{r^{2h(t)}}\int_{\bbbr} \int_{\frac{K}{r}} \left| \frac{C_{\alpha(t+r)}^{1/ \alpha(t+r)}}{y^{1/ \alpha(t+r)}} f(t+r,t+r,x) - \frac{C_{\alpha(t)}^{1/ \alpha(t)}}{y^{1/ \alpha(t)}} f(t+r,t,x) \right|^2 dy \hspace*{0.1cm} m(dx) ,
\end{displaymath}
\begin{displaymath}
\Delta_2(r,t) = \frac{2C_{\alpha(t)}^{1/ \alpha(t)}}{r^{2h(t)}}  \int_{\bbbr} \int_{\frac{K}{r}} \frac{1}{y^{1/ \alpha(t)}} g_1(r,t,x,y) g_2(r,t,x) dy \hspace*{0.1cm} m(dx),
\end{displaymath}
and
\begin{displaymath}
\Delta_3(r,t) = \frac{1}{r^{2h(t)}}\int_{\bbbr} \int_{\frac{K}{r}} \frac{C_{\alpha(t)}^{2/ \alpha(t)}}{y^{2/ \alpha(t)}} \left(f(t+r,t,x) - f(t,t,x) \right)^2 dy \hspace*{0.1cm} m(dx) - l(t),
\end{displaymath}
with $g_1(r,t,x,y) = \frac{C_{\alpha(t+r)}^{1/ \alpha(t+r)}}{y^{1/ \alpha(t+r)}} f(t+r,t+r,x) - \frac{C_{\alpha(t)}^{1/ \alpha(t)}}{y^{1/ \alpha(t)}} f(t+r,t,x)$ and $ g_2(r,t,x) = f(t+r,t,x) - f(t,t,x)$. Since $\alpha$ is continuous, there exists a positive constant $K_U$ (that may change from line to line) such that 

\begin{eqnarray*}
|\Delta_2(r,t) | & \leq & \frac{K_U}{r^{2h(t)}} \int_{\bbbr} \int_{\frac{K}{r}} \left| \frac{g_1(r,t,x,y) g_2(r,t,x)}{y^{1/ \alpha(t)}}\right| dy \hspace*{0.1cm} m(dx)\\
& \leq & \frac{K_U}{r^{2h(t)}} \left( \int_{\bbbr} \int_{\frac{K}{r}}  \left| g_1(r,t,x,y) \right|^2 dy \hspace*{0.1cm} m(dx)\right)^{\frac{1}{2}} \left( \int_{\bbbr} \int_{\frac{K}{r}}  \left| \frac{g_2(r,t,x)}{y^{1/ \alpha(t)}} \right|^2 dy \hspace*{0.1cm} m(dx) \right)^{\frac{1}{2}} \\
& \leq & \frac{K_U}{r^{2h(t)}} r^{h(t)} \sqrt{\Delta_1(r,t)} \left( \int_{\bbbr} \int_{\frac{K}{r}}  \left| \frac{g_2(r,t,x)}{y^{1/ \alpha(t)}} \right|^2 dy \hspace*{0.1cm} m(dx) \right)^{\frac{1}{2}} \\
& \leq & \frac{K_U}{r^{h(t)}} \sqrt{\Delta_1(r,t)} \left( \int_{\bbbr} \left| g_2(r,t,x) \right|^2  m(dx) \right)^{\frac{1}{2}} r^{\frac{1}{\alpha(t)} - \frac{1}{2}} K^{ \frac{1}{2} -\frac{1}{\alpha(t)}} \sqrt{\frac{\alpha(t)}{2-\alpha(t)}}\\
& \leq & K_U \sqrt{\Delta_1(r,t)} \left( \frac{1}{r^{1+2(h(t)-\frac{1}{\alpha(t)})}} \int_{\bbbr} \left| f(t+r,t,x) - f(t,t,x) \right|^2 m(dx) \right)^{\frac{1}{2}} \\
& \leq &  K_U \sqrt{\Delta_1(r,t)} \textrm{ with } \textrm{\bf  (\ref{C16})}.\\
\end{eqnarray*}
Let us show that $ \lim_{r \rightarrow 0}  \sqrt{\Delta_1(r,t)} =0.$ The triangle inequality yields $\sqrt{\Delta_1(r,t)}) \leq \delta_1(r,t) + \delta_2(r,t) + \delta_3(r,t)$ where 
\begin{displaymath}
\delta_1(r,t) = \frac{1}{2r^{h(t)}} \left( \int_{\bbbr} \int_{\frac{K}{r}} \left| C_{\alpha(t+r)}^{1/ \alpha(t+r)} - C_{\alpha(t)}^{1/ \alpha(t)} \right|^2 \frac{|f(t+r,t+r,x)|^2}{y^{2/ \alpha(t+r)}} dy \hspace*{0.1cm} m(dx) \right)^{\frac{1}{2}} ,
\end{displaymath}
\begin{displaymath}
\delta_2(r,t) = \frac{1}{2r^{h(t)}} \left( \int_{\bbbr} \int_{\frac{K}{r}} \frac{C_{\alpha(t)}^{2/ \alpha(t)}}{y^{2/ \alpha(t+r)}} \left| f(t+r,t+r,x) - f(t+r,t,x)\right|^2  dy \hspace*{0.1cm} m(dx) \right)^{\frac{1}{2}} ,
\end{displaymath}
and
\begin{displaymath}
\delta_3(r,t) = \frac{1}{2r^{h(t)}} \left( \int_{\bbbr} \int_{\frac{K}{r}} C_{\alpha(t)}^{2/ \alpha(t)}|f(t+r,t,x)|^2 \left( \frac{1}{y^{1/ \alpha(t+r)}} - \frac{1}{y^{1/ \alpha(t)}}\right)^2 dy \hspace*{0.1cm} m(dx) \right)^{\frac{1}{2}}.
\end{displaymath}
Now, 
\begin{eqnarray*}
\delta_1(r,t) & \leq & K_U \frac{| C_{\alpha(t+r)}^{1/ \alpha(t+r)} - C_{\alpha(t)}^{1/ \alpha(t)} | }{r^{h(t)}} \left( \frac{1}{\frac{2}{\alpha(t+r)}-1} (\frac{K}{r} )^{1-\frac{2}{\alpha(t+r)}} \right)^{\frac{1}{2}} \left( \int_{\bbbr}  \left| f(t+r,t+r,x) \right|^2 m(dx) \right)^{\frac{1}{2}}. \\
\end{eqnarray*}
Since the function $u \mapsto C^{1/\alpha(u)}_{\alpha(u)}$ is a $C^1$ function, 
\begin{eqnarray*}
\delta_1(r,t) & \leq & K_U r^{1-h(t)+\frac{1}{\alpha(t+r)}-\frac{1}{2}}  \left( \int_{\bbbr}  \left| f(t+r,t+r,x) \right|^2 m(dx) \right)^{\frac{1}{2}} \\
& \leq & K_U r^{1-h(t)+\frac{1}{\alpha(t+r)}-\frac{1}{2}} \textrm{ with } \textrm{\bf  (C12)}\\
& \leq & K_U r^{\frac{1}{2} + \frac{1}{d} - h_+}.\\
\end{eqnarray*}
Since $\frac{1}{2} + \frac{1}{d} - h_+ > 0 $, $\lim_{r \rightarrow 0}\limits  \delta_1(r,t) \ = 0$.
\begin{eqnarray*}
\delta_2(r,t) & \leq & \frac{C_{\alpha(t)}^{1/ \alpha(t)}}{2r^{h(t)}} \left( \frac{1}{\frac{2}{\alpha(t+r)}-1} (\frac{K}{r} )^{1-\frac{2}{\alpha(t+r)}} \right)^{\frac{1}{2}}  \left( \int_{\bbbr}  \left| f(t+r,t+r,x) - f(t+r,t,x) \right|^2 m(dx) \right)^{\frac{1}{2}} \\
& \leq & K_U r^{\frac{1}{\alpha(t+r)}-h(t)-\frac{1}{2}} \left( \int_{\bbbr}  \left| f(t+r,t+r,x) - f(t+r,t,x) \right|^2 m(dx) \right)^{\frac{1}{2}} \\
& \leq & K_U r^{\frac{1}{2}+ \frac{1}{\alpha(t+r)}-h(t)} \textrm{ with } \textrm{\bf  (C15)} \\
& \leq & K_U r^{\frac{1}{2} + \frac{1}{d} - h_+},\\
\end{eqnarray*}
thus $\lim_{r \rightarrow 0}\limits  \delta_2(r,t)  = 0$.

\begin{eqnarray*}
\delta_3(r,t) & \leq & \frac{C_{\alpha(t)}^{1/ \alpha(t)}}{2r^{h(t)}}  \left( \int_{\bbbr}  \left| f(t+r,t,x) \right|^2 m(dx) \right)^{\frac{1}{2}} \left( \int_{\frac{K}{r}} \left( \frac{1}{y^{1/ \alpha(t+r)}} - \frac{1}{y^{1/ \alpha(t)}}\right)^2 dy \right)^{\frac{1}{2}} \\
\end{eqnarray*}

Since the function $u \mapsto \alpha(u)$ is a $C^1$ function, $\forall \eta < \frac{1}{d}$, 
\begin{eqnarray*}
\delta_3(r,t) & \leq & \frac{K_U}{r^{h(t)}} \left( \int_{\bbbr}  \left| f(t+r,t,x) \right|^2 m(dx) \right)^{\frac{1}{2}} K_U r^{\frac{1}{2} + \frac{1}{d}- \eta }\\
& \leq & K_U r^{\frac{1}{2} + \frac{1}{d}- \eta - h_+} \textrm{ with } \textrm{\bf  (C11)}\\
\end{eqnarray*}
thus $\lim_{r \rightarrow 0}\limits  \delta_3(r,t)  = 0$. Finally, $\lim_{r \rightarrow 0}\limits \sqrt{\Delta_1(r,t)}  =0.$ 

Let us now consider the last term $\Delta_3(r,t)$:
\begin{eqnarray*}
\Delta_3(r,t) & = & \frac{C_{\alpha(t)}^{2/ \alpha(t)} K^{1- \frac{2}{\alpha(t)}}}{\frac{2}{\alpha(t)}-1} \left( \frac{1}{r^{1+2(h(t)-1/\alpha(t)}} \int_{\bbbr} \left( f(t+r,t,x) - f(t,t,x) \right)^2  m(dx) - g(t) \right)\\
\end{eqnarray*}
thus, with {\bf  (C14)}, $\lim_{r \rightarrow 0}\limits \left| \Delta_3(r,t) \right| =0$ \Box

\begin{lem} \label{lemvit}
Assume (C6), (C10), (C12), (C13), (C15), and let:
\begin{displaymath}
\Delta(r,t) =: \frac{1}{r^{1+2(h(t)-1/ \alpha(t))}} \left( \frac{C_{\alpha(t)}^{1/ \alpha(t)} K^{\frac{1}{ \alpha(t+r)} - \frac{1}{\alpha(t)}}  (\frac{2}{\alpha(t+r)} -1 ) \int_{\bbbr} f(t+r,t+r,x)f(t,t,x) m(dx)}{C_{\alpha(t+r)}^{1/ \alpha(t+r)} r^{\frac{1}{ \alpha(t+r)} - \frac{1}{ \alpha(t)}} ( \frac{1}{\alpha(t+r)} + \frac{1}{\alpha(t)} -1) \int_{\bbbr} f(t+r,t+r,x)^2 m(dx) } -1 \right)^2.
\end{displaymath}
Then:
\begin{displaymath}
\lim_{r \rightarrow 0} |\Delta(r,t)| =0.
\end{displaymath}

If in addition  we suppose (Cu10), (Cu12), (Cu15), the convergence is uniform on $U$.
\end{lem}
 
\noindent
{\bf Proof} 

\noindent
Since the function $t \mapsto \alpha(t)$ is a $C^1$ function, there exists $K_U > 0$ such that 
 \begin{equation}\label{ineq1}
 \left| \frac{C_{\alpha(t)}^{1/ \alpha(t)}}{C_{\alpha(t+r)}^{1/ \alpha(t+r)}} - 1\right| \leq r K_U,
 \end{equation} 
  \begin{equation}\label{ineq2} 
 \left| K^{\frac{1}{ \alpha(t+r)} - \frac{1}{\alpha(t)}} - 1\right| \leq r K_U,
 \end{equation}
 and
  \begin{equation}\label{ineq3}
 \left| \frac{\frac{2}{\alpha(t+r)} -1}{\frac{1}{\alpha(t+r)} + \frac{1}{\alpha(t)} -1} - 1\right| \leq r K_U.
 \end{equation} 
 Increasing $K_U$ if necessary, we also have, $\forall a >0$, 
  \begin{equation}\label{ineq4}
 \left| \frac{1}{r^{\frac{1}{\alpha(t+r)}- \frac{1}{\alpha(t)}}} - 1\right| \leq r^a K_U.
 \end{equation} 
 For the last term, we write
 \begin{displaymath}
 \frac{\int_{\bbbr} f(t+r,t+r,x)f(t,t,x) m(dx)}{\int_{\bbbr} f(t+r,t+r,x)^2 m(dx)} - 1 = \Delta_1(r,t) + \Delta_2(r,t)
 \end{displaymath}
 where
 \begin{displaymath}
 \Delta_1(r,t) = \frac{1}{\int_{\bbbr} f(t+r,t+r,x)^2 m(dx)} \left( \int_{\bbbr} f(t+r,t+r,x) \left( f(t,t,x) - f(t+r,t,x)\right) m(dx)\right)
 \end{displaymath}
 and
 \begin{displaymath}
 \Delta_2(r,t) = \frac{1}{\int_{\bbbr} f(t+r,t+r,x)^2 m(dx)} \left( \int_{\bbbr}  f(t+r,t+r,x)  \left( f(t+r,t,x) - f(t+r,t+r,x)\right) m(dx) \right).
 \end{displaymath}
 With {\bf (C13)}, we may choose $K_U$ such that 
 \begin{displaymath}
 |\Delta_1(r,t)| \leq K_U \int_{\bbbr} |f(t+r,t+r,x)| \left| f(t,t,x) - f(t+r,t,x)\right| m(dx) .
 \end{displaymath}
 Let $p \in (\alpha(t),2)$ , $ p \geq 1$ satisfying {\bf (C10)}, and $q $ such that $\frac{1}{p} + \frac{1}{q} =1$. H\"older inequality entails:
 \begin{eqnarray*}
  |\Delta_1(r,t)| & \leq & K_U \left( \int_{\bbbr} |f(t+r,t+r,x)|^q m(dx)\right)^{1/ q } \left(\int_{\bbbr}  \left| f(t,t,x) - f(t+r,t,x)\right|^p  m(dx)\right)^{1/ p}\\
  & \leq & K_U \left(\int_{\bbbr}  \left| f(t+r,t,x) - f(t,t,x)\right|^p  m(dx)\right)^{1/ p} \textrm{ with } \textrm{\bf  (C6)} \textrm{ and } \textrm{\bf (C12)}\\
  & \leq & K_U r^{\frac{1}{p}+h(t)-\frac{1}{\alpha(t)}} \textrm{ with } \textrm{\bf  (C10)}. \\
 \end{eqnarray*} 
With {\bf (C12)}, {\bf (C13)} and Cauchy-Schwarz inequality, we select $K_U$ such that 
\begin{eqnarray*}
 |\Delta_2(r,t)| & \leq & K_U \left( \int_{\bbbr} \left| f(t+r,t+r,x) - f(t+r,t,x)\right|^2  m(dx) \right)^{\frac{1}{2}}\\
 & \leq & K_U r \textrm{ with } \textrm{\bf  (C15)}.\\
\end{eqnarray*}
Finally, since $ h(t) +\frac{1}{p} - \frac{1}{\alpha(t)} \leq 1$,
\begin{equation}\label{ineq5}
\left| \frac{\int_{\bbbr} f(t+r,t+r,x)f(t,t,x) m(dx)}{\int_{\bbbr} f(t+r,t+r,x)^2 m(dx)} - 1 \right| \leq K_U r^{h(t) +\frac{1}{p} - \frac{1}{\alpha(t)}}.
\end{equation}
Using the inequalities (\ref{ineq1}), (\ref{ineq2}), (\ref{ineq3}), (\ref{ineq4}) and (\ref{ineq5}), we may find a constant $K_U$ such that for all $a>0$,
\begin{displaymath}
 |\Delta(r,t)| \leq \frac{1}{r^{1+2(h(t)-1/ \alpha(t))}} K_U (r^2 + r^{2a} + r^{2(h(t) +\frac{1}{p} - \frac{1}{\alpha(t)})}).
\end{displaymath}

Choosing $a \in \left( h(t) +\frac{1}{p} - \frac{1}{\alpha(t)} , 1\right)$, this entails:
\begin{eqnarray*}
|\Delta(r,t)| & \leq &  \frac{3}{r^{1+2(h(t)-1/ \alpha(t))}} K_U r^{2(h(t) +\frac{1}{p} - \frac{1}{\alpha(t)})}\\
& \leq & 3 K_U r^{\frac{2}{p}-1}.\\
\end{eqnarray*}
Since $\frac{2}{p}-1 > 0$, $\lim_{r \rightarrow 0}\limits |\Delta(r,t)| =0 $ \Box

\begin{lem}\label{lemsin}
Assuming (C1), (C6), (C7), (C8), one has:
$$\forall \varep < \frac{1}{d},  \exists K_U \leq 1 \mbox{ such that }\forall v \geq 1 ,  \forall r \leq \varep_0,$$
\begin{eqnarray*}
y \geq K_U\frac{ v^{\frac{d}{1-\varep d}}}{r} & \Rightarrow & \sin^2 \left( \frac{v C_{\alpha(t+r)}^{1/ \alpha(t+r)}f(t+r,t+r,x)}{2r^{h(t)}y^{1/ \alpha(t+r)}}  - \frac{v C_{\alpha(t)}^{1/ \alpha(t)}f(t,t,x)}{2r^{h(t)}y^{1/ \alpha(t)}} \right) \\
& \geq & \frac{1}{2} \left| \frac{v C_{\alpha(t+r)}^{1/ \alpha(t+r)}f(t+r,t+r,x)}{2r^{h(t)}y^{1/ \alpha(t+r)}}  - \frac{v C_{\alpha(t)}^{1/ \alpha(t)}f(t,t,x)}{2r^{h(t)}y^{1/ \alpha(t)}} \right|^2.\\
\end{eqnarray*}
If in addition we suppose (Cu8),
\begin{eqnarray*}
y \geq K_U\frac{ v^{\frac{d}{1-\varep d}}}{r} & \Rightarrow & \forall t \in U, \sin^2 \left( \frac{v C_{\alpha(t+r)}^{1/ \alpha(t+r)}f(t+r,t+r,x)}{2r^{h(t)}y^{1/ \alpha(t+r)}}  - \frac{v C_{\alpha(t)}^{1/ \alpha(t)}f(t,t,x)}{2r^{h(t)}y^{1/ \alpha(t)}} \right) \\
& \geq & \frac{1}{2} \left| \frac{v C_{\alpha(t+r)}^{1/ \alpha(t+r)}f(t+r,t+r,x)}{2r^{h(t)}y^{1/ \alpha(t+r)}}  - \frac{v C_{\alpha(t)}^{1/ \alpha(t)}f(t,t,x)}{2r^{h(t)}y^{1/ \alpha(t)}} \right|^2.\\
\end{eqnarray*}
\end{lem}

\noindent
{\bf Proof} 

\noindent
Let $\varep < \frac{1}{d}$. We write $ \frac{v C_{\alpha(t+r)}^{1/ \alpha(t+r)}f(t+r,t+r,x)}{2r^{h(t)} y^{1/ \alpha(t+r)}}  - \frac{v C_{\alpha(t)}^{1/ \alpha(t)} f(t,t,x)}{2r^{h(t)} y^{1/ \alpha(t)}} = \kappa_1 (r,t,v,x,y) + \kappa_2 (r,t,v,x,y)$, with
\begin{displaymath}
\kappa_1(r,t,v,x,y) = \frac{v}{2r^{h(t)}} \left( \frac{ C_{\alpha(t+r)}^{1/ \alpha(t+r)}f(t+r,t+r,x)}{ y^{1/ \alpha(t+r)}}  - \frac{ C_{\alpha(t)}^{1/ \alpha(t)} f(t+r,t,x)}{ y^{1/ \alpha(t)}} \right)
\end{displaymath}
and 
\begin{displaymath}
\kappa_2(r,t,v,x,y) = \frac{v C_{\alpha(t)}^{1/ \alpha(t)}}{2r^{h(t)} y^{1/ \alpha(t)}} \left( f(t+r,t,x) - f(t,t,x) \right).
\end{displaymath}
Using the finite-increments theorem,
\begin{eqnarray*}
\left| \kappa_1(r,t,v,x,y) \right| & \leq & \frac{v}{2r^{h(t)}} r  ( \sup_{a \in U}\left| \frac{K_U | f(t+r,a,x)|}{y^{1/ \alpha(a)}} \right| + \sup_{a \in U}\left| \frac{C_{\alpha(a)}^{1/ \alpha(a)} | f_v(t+r,a,x)|}{y^{1/ \alpha(a)}}\right| \\
& + & \sup_{a \in U}\left| \frac{| \alpha'(a)|}{\alpha^2(a)} |\ln y| \frac{C_{\alpha(a)}^{1/ \alpha(a)} | f(t+r,a,x)|}{y^{1/ \alpha(a)}} \right| ) .
\end{eqnarray*}
For $y \geq 1$, conditions {\bf (C6)} and {\bf (C7)} imply
\begin{displaymath}
\frac{K_U | f(t+r,a,x)|}{y^{1/ \alpha(a)}} \leq \frac{K_U}{y^{1/ d}},
\end{displaymath}

\begin{displaymath}
\frac{K_U | f_v(t+r,a,x)|}{y^{1/ \alpha(a)}} \leq \frac{K_U}{y^{1/ d}},
\end{displaymath}
and
\begin{displaymath}
 \frac{| \alpha'(a)|}{\alpha^2(a)} |\ln y| \frac{C_{\alpha(a)}^{1/ \alpha(a)} | f(t+r,a,x)|}{y^{1/ \alpha(a)}} \leq \frac{K_U |\ln y|}{y^{1/ d}}.
\end{displaymath}

Finally,
\begin{eqnarray*}
\left| \kappa_1(r,t,v,x,y) \right| & \leq &  \frac{K_U v r^{1-h(t)}}{y^{1/ d}} \left( 1 + |\ln y|\right)\\
& \leq & \frac{K_U v}{y^{1/d - \varep}}. \\
\end{eqnarray*}
Condition {\bf (C8)} allows to estimate $\kappa_2(r,t,v,x,y)$ as follows: 
\begin{displaymath}
\left| \kappa_2(r,t,v,x,y)\right| \leq \frac{K_U v}{(ry)^{1/ \alpha(t)}}.
\end{displaymath}
Finally, $\forall K >0$, $\forall \varep < \frac{1}{d}$,$\exists K_U \geq 1$,$\forall v \geq 1 $,$\forall r < \varep_0$, $\forall y \geq K_U \frac{v^{\frac{d}{1-\varep d}}}{r}$,
\begin{displaymath}
\left| \frac{v C_{\alpha(t+r)}^{1/ \alpha(t+r)}f(t+r,t+r,x)}{2r^{h(t)} y^{1/ \alpha(t+r)}}  - \frac{v C_{\alpha(t)}^{1/ \alpha(t)} f(t,t,x)}{2r^{h(t)} y^{1/ \alpha(t)}} \right| \leq K
\end{displaymath}
\Box

\begin{lem}\label{borneN}
Assuming  (C6), (C10), (C11), (C12), (C13), (C14), (C15), there exist $\varep_0$ and $K_U >0$ such that $\forall r < \varep_0$, $\forall v \geq 1$:
\begin{displaymath}
N(v,t,r) \geq K_U v^{2+\frac{d}{1-\varep d}(1- \frac{2}{c})},
\end{displaymath}
where
\begin{displaymath}
 N(v,t,r) = : \int_{\bbbr} \int_{\frac{K_U v^{\frac{d}{1-\varep d}}}{r}} \left| \frac{v C_{\alpha(t+r)}^{1/ \alpha(t+r)}f(t+r,t+r,x)}{2r^{h(t)}y^{1/ \alpha(t+r)}}  - \frac{v C_{\alpha(t)}^{1/ \alpha(t)}f(t,t,x)}{2r^{h(t)}y^{1/ \alpha(t)}} \right|^2 dy \hspace*{0.1cm} m(dx).
\end{displaymath}
If in addition we suppose (Cu10), (Cu11), (Cu12), (Cu14), (Cu15), the constant $K_U$ does not depend on $t$.
\end{lem}

\noindent
{\bf Proof} 

\noindent
Expanding the square above, we may write 
\begin{displaymath}
 N(v,t,r) = A_1(r,t) v^{2+\frac{d}{1-\varep d}(1-\frac{2}{\alpha(t+r)})} - A_2(r,t) v^{2+\frac{d}{1-\varep d}(1-\frac{1}{\alpha(t+r)}-\frac{1}{\alpha(t)})} + A_3(r,t) v^{2+\frac{d}{1-\varep d}(1-\frac{2}{\alpha(t)})},
\end{displaymath}
where
\begin{displaymath}
 A_1(r,t) = \frac{C_{\alpha(t+r)}^{2/ \alpha(t+r)}  (K_U)^{1-\frac{2}{\alpha(t+r)}}}{4\left( \frac{2}{\alpha(t+r)} -1\right) r^{1+2(h(t)-\frac{1}{\alpha(t+r)})}} \int_{\bbbr} |f(t+r,t+r,x)|^2 m(dx),
\end{displaymath}
\begin{displaymath}
 A_2(r,t) = \frac{C_{\alpha(t+r)}^{1/ \alpha(t+r)} C_{\alpha(t)}^{1/ \alpha(t)}  (K_U)^{1-\frac{1}{\alpha(t+r)}-\frac{1}{\alpha(t)}}}{2\left( \frac{1}{\alpha(t+r)} +\frac{1}{\alpha(t)} -1\right) r^{1+2h(t)-\frac{1}{\alpha(t+r)} - \frac{1}{\alpha(t)}}} \int_{\bbbr} f(t+r,t+r,x) f(t,t,x) m(dx),
\end{displaymath}
and
\begin{displaymath}
 A_3(r,t) = \frac{C_{\alpha(t)}^{2/ \alpha(t)}  (K_U)^{1-\frac{2}{\alpha(t)}}}{4\left( \frac{2}{\alpha(t)} -1\right) r^{1+2(h(t)-\frac{1}{\alpha(t)})}} \int_{\bbbr} |f(t,t,x)|^2 m(dx).
\end{displaymath}
We obtain 
\begin{displaymath}
N(v,t,r) = v^{2+\frac{d}{1-\varep d}(1- \frac{2}{\alpha(t)})} \left( A_1(r,t) (v^{\frac{2d}{1-\varep d}(\frac{1}{\alpha(t)}-\frac{1}{\alpha(t+r)})})^2 - A_2(r,t)(v^{\frac{2d}{1-\varep d}(\frac{1}{\alpha(t)}-\frac{1}{\alpha(t+r)})})+ A_3(r,t) \right).
\end{displaymath}
Let $P(r,t,X)=A_1(r,t) X^2 - A_2(r,t)X+ A_3(r,t)$. Write: 
$$P(r,t,X) = P(r,t,X) - P(r,t,\frac{A_2(r,t)}{2A_1(r,t)}) + P(r,t,\frac{A_2(r,t)}{2A_1(r,t)}) -P(r,t,1)+P(r,t,1).$$
Since $P(\frac{A_2(r,t)}{2A_1(r,t)}) $ is the minimum of $P$, 
\begin{displaymath}
P(r,t,X) \geq P(r,t,\frac{A_2(r,t)}{2A_1(r,t)}) -P(r,t,1)+P(r,t,1).
\end{displaymath}
Note that $P(r,t,1) =N(1,t,r)$, thus lemma (\ref{lemN}) entails that there exists a positive function $l$ such that $\lim_{r \rightarrow 0}\limits P(r,t,1) = l(t)$. For $P(r,t,\frac{A_2(r,t)}{2A_1(r,t)}) -P(r,t,1)$, we use lemma (\ref{lemvit}). With the same notations,
\begin{eqnarray*}
|P(r,t,\frac{A_2(r,t)}{2A_1(r,t)}) -P(r,t,1)| & = & | A_1(r,t)| r^{1+2(h(t)-\frac{1}{\alpha(t)})} \Delta(r,t)\\
& \leq & K_U \Delta(r,t)\\
\end{eqnarray*}
thus $\lim_{r \rightarrow 0}\limits |P(r,t,\frac{A_2(r,t)}{2A_1(r,t)}) -P(r,t,1)| =0$. As a
consequence, there exist a positive constant $K_U$ and $\varep_0 >0$ such that for all $x \in \bbbr$ and $r \in (0, \varep_0)$, $P(r,t,x) \geq K_U$. We obtain $N(v,t,r) \geq v^{2+\frac{d}{1-\varep d}(1- \frac{2}{\alpha(t)})} K_U$ for all $v \in \bbbr$ and $r \in (0, \varep_0)$. Since $\alpha(t) >c$, $N(v,t,r) \geq K_U v^{2+\frac{d}{1-\varep d}(1- \frac{2}{c})}$ \Box
\newline

\noindent
{\bf Proof of theorem \ref{vitesp}}

\noindent
Consider
\begin{displaymath}
 \E \left[ \left|\frac{Y(t+\varep)-Y(t)}{\varep^{h(t)}}\right|^{\eta}\right] = \int_{0}^{\infty} \P \left( \left|\frac{Y(t+\varep)-Y(t)}{\varep^{h(t)}}\right|^{\eta} > x \right)dx.
\end{displaymath}
Thanks to {\bf  (C1)}, {\bf  (C2)}, {\bf  (C3)} and {\bf  (C5)}, $Y$ is $h(t)$-localisable at $t$  \cite{LGLV}, thus for all $x>0$,
\begin{displaymath}
 \P \left( \left|\frac{Y(t+\varep)-Y(t)}{\varep^{h(t)}}\right|^{\eta} > x \right) \rightarrow \P \left( \left|Y'_t(1)\right|^{\eta} > x \right).
\end{displaymath}
We shall make use of Lebesgue dominated convergence theorem.
\newline

For $x \leq 1$, $\P \left( \left|\frac{Y(t+\varep)-Y(t)}{\varep^{h(t)}}\right|^{\eta} > x \right) \leq 1$.
\newline

For $x >1$,
\begin{eqnarray*}
 \P \left( \left|\frac{Y(t+\varep)-Y(t)}{\varep^{h(t)}}\right|^{\eta} > x \right) &=& \P \left( \left|\frac{Y(t+\varep)-Y(t)}{\varep^{h(t)}}\right| > x^{1/ \eta} \right)\\
 & \leq & \P \left(\left| \frac{X(t+\varep,t+\varep)-X(t+\varep,t)}{\varep^{h(t)}}\right|>\frac{x^{1/\eta}}{2} \right)\\
 & & + \hspace*{0.3cm} \P \left(\left|\frac{X(t+\varep,t)-X(t,t)}{\varep^{h(t)}} \right| > \frac{x^{1/ \eta}}{2} \right).\\
\end{eqnarray*}
For the first term, by proposition \ref{msspfm} (or \ref{msspfm2}),
\begin{eqnarray*}
\P \left(\left| \frac{X(t+\varep,t+\varep)-X(t+\varep,t)}{\varep^{h(t)}}\right|>\frac{x^{1/\eta} }{2} \right) & \leq & \frac{K_U}{x^{d/\eta}} \left( 1 + |\log x|^d \right) + \frac{K_U}{x^{c/\eta}} \left( 1 + |\log x|^c \right).\\
\end{eqnarray*}
For the second term, let $p \in (\eta,\alpha(t))$.

\begin{eqnarray*}
 \P \left(\left|\frac{X(t+\varep,t)-X(t,t)}{\varep^{h(t)}} \right| > \frac{x^{1/ \eta}}{2} \right) & = & \P \left( \left|\frac{X(t+\varep,t)-X(t,t)}{\varep^{h(t)}} \right|^p > \frac{x^{p/\eta}}{2^p} \right).\\
\end{eqnarray*}
With Markov inequality and {\bf  (C9)},
\begin{eqnarray*}
 \P \left(\left|\frac{X(t+\varep,t)-X(t,t)}{\varep^{h(t)}} \right| > \frac{x^{1/ \eta}}{2} \right ) & \leq & \frac{2^p}{x^{p/\eta}\varep^{ph(t)}} C_{\alpha(t),0}(p)^p \Vert f(t+\varep,t,.)-f(t,t,.) \Vert_{\alpha(t)}^p\\
 & \leq & \frac{2^p C_{\alpha(t),0}(p)^p}{x^{p/\eta} \varep^{ph(t)}} \left(\int_{\bbbr} | f(t+\varep,t,x)-f(t,t,x)|^{\alpha(t)} m(dx) \right)^{p/ \alpha(t)}\\
 & \leq & \frac{K_{p,\alpha(t)}}{x^{p/\eta}},\\
 \end{eqnarray*}
 
 thus
\begin{displaymath}
 \P( \left|\frac{Y(t+\varep)-Y(t)}{\varep^{h}}\right|^{\eta} > x) \leq K_U \left(\frac{1}{x^{d/\eta}} \left( 1 + |\log x|^d \right) + \frac{1}{x^{c/\eta}} \left( 1 + |\log x|^c \right)+\frac{1}{x^{p/\eta}} \right) \mathbf{1}_{x>1}+\mathbf{1}_{x \leq 1} \Box
\end{displaymath}



\noindent
{\bf Proof of theorem \ref{upbound}} 

\noindent
Let $ \gamma  > h(t)$ and $x >0$. 
\begin{displaymath}
\P \left( \frac{r^{\gamma}}{|Y(t+r) -Y(t)|} > x \right) = \P \left( | Y(t+r) -Y(t)| < \frac{r^{\gamma}}{x}\right).
\end{displaymath}
Applying proposition (\ref{intfonc}), there exists $\varep_0 > 0 $ such that 
\begin{displaymath}
\sup_{r \in B(0,\varep_0)} \int_{0}^{+\infty} \varphi_{\frac{Y(t+r)-Y(t)}{r^{h(t)}}}(v) dv <  + \infty.
\end{displaymath}

Thus with proposition (\ref{probinf}), there exists $K_U > 0$ such that 
\begin{displaymath}
 \P \left( | Y(t+r) -Y(t)| < \frac{r^{\gamma}}{x}\right) \leq K_U \frac{r^{\gamma - h(t)}}{x}.
\end{displaymath}

Let $r_n = \frac{1}{n^{\eta}}$ with $\eta (\gamma - h(t)) >1$. $\forall x >0$, $\sum_{n}\limits \P \left( \frac{r_n^{\gamma}}{|Y(t+r_n) -Y(t)|} > x \right) < +\infty$. Borel Cantelli lemma entails
that, almost surely, $\lim_{n \rightarrow +\infty}\limits \frac{|Y(t+r_n) -Y(t)|}{r_n^{\gamma}} = +\infty$. As a consequence, almost surely, $\limsup_{r \rightarrow 0}\limits \frac{|Y(t+r) -Y(t)|}{r^{\gamma}} = +\infty$, and 
\begin{displaymath}
\h \leq h(t). 
\end{displaymath}
\Box

\noindent
{\bf Proof of theorem \ref{esplevy}} 

\noindent
We want to apply theorem (\ref{vitesp}) with $f(t,u,x) = \mathbf{1}_{[0,t]}(x)$. Let us show that conditions (C1), (C2), (C3), (C5) and (C9) are satisfied.
\begin{itemize}
    \item (C1) The family of functions $v \to f(t,v,x)$ is differentiable for all $(v,t)$ in $(0,1)^2$ and almost all $x$ in $E$. The derivatives of $f$ with respect to $u$ vanish.
    \item (C2) 
    \begin{displaymath}
  |f(t,w,x)|^{\alpha(w)} = \mathbf{1}_{[0,t]}(x)   
    \end{displaymath}
thus, for all $\delta > 0$, all $t \in (0,1)$,
\begin{displaymath}
\int_\bbbr \left[\sup_{w \in (0,1)} (|f(t,w,x)|^{\alpha(w)})\right]^{1+\delta} \hspace{0.1cm} dx = t
\end{displaymath}
and (C2) holds.

\item (C3) $f'_u = 0$ thus (C3) holds.

\item (C5) $X(t,u)$ (as a process in $t$) is localisable at $u$ with exponent $\frac{1}{\alpha(u)} \in (\frac{1}{d},\frac{1}{c}) \subset (0,1)$, with local form $X_u(t,u)$, and  $u \mapsto \frac{1}{\alpha(u)}$ is a $C^1$ function (see \cite{LGLV}).

\item (C9) 
\begin{eqnarray*}
 \frac{1}{r^{h(t)\alpha(t)}} \int_{\bbbr} | f(t+r,t,x)-f(t,t,x)|^{\alpha(t)} m(dx) & = & \frac{1}{r} \int_{t}^{t+r} dx \\
 & = & 1,\\
\end{eqnarray*}
thus (C9) holds.
\end{itemize}

From theorem (\ref{vitesp}), we get that 
 \begin{displaymath}
 \E \left[ |Y(t+\varep)-Y(t)|^{\eta}\right] \sim \varep^{\frac{\eta}{\alpha(t)}} \E \left[ |Y'_t(1)|^{\eta}\right].
\end{displaymath}
Since $Y'_t(1)$ is an $S_{\alpha(t)}(1,0,0)$ random variable, property 1.2.17 of \cite{ST} allows to conclude. \Box

\noindent
{\bf Proof of theorem \ref{expolevy}} 

\noindent
We want to apply Theorem (\ref{upbound}) with $f(t,u,x) = \mathbf{1}_{[0,t]}(x)$ and $h(t)= \frac{1}{\alpha(t)}$ in order to obtain the inequality. Let us show that the conditions (C6), (C7), (Cu8), (Cu10), (Cu11), (Cu12), (C13), (Cu14) and (Cu15) are satisfied.
\begin{itemize}
    \item (C6) Obvious.
    \item (C7) Obvious.
    \item (Cu8) $\forall v \in U$, $\forall u \in U$, $\forall x \in \bbbr$,
    \begin{eqnarray*}
    \frac{1}{|v-u|^{h(u)-1/ \alpha(u)}}\left| f(v,u,x) -f(u,u,x) \right|  & = & \mathbf{1}_{[u,v]}(x)\\
    & \leq & 1\\
    \end{eqnarray*}
    thus (Cu8) holds.
    \item (Cu10) $\forall v \in U$,$\forall u \in U$,
    \begin{eqnarray*}
     \frac{1}{|v-u|^{1+p(h(u)-\frac{1}{\alpha(u)})}} \int_{\bbbr} \left| f(v,u,x) - f(u,u,x) \right|^p m(dx) & = & \frac{1}{|v-u|} \int_{\bbbr} | \mathbf{1}_{[u,v]} (x)|\\
 & = & 1\\
 \end{eqnarray*}
 thus (Cu10) holds.
    \item (Cu11) $\forall  v \in U$, $\forall  u \in U$, 
	\begin{displaymath}
	\int_{\bbbr} \left| f(v,u,x) \right|^2 m(dx) = v
	\end{displaymath}
	thus (Cu11) holds ($U=(0,1)$).
    \item (Cu12) For the same reason as (Cu11), (Cu12) holds.
    \item (C13) Since $t \in (0,1)$ (in particular $ t \neq 0$), one can choose $U$ such that $\inf_{v \in U} v > 0$ thus (C13) holds.
    \item (Cu14) 
    \begin{eqnarray*}
    \frac{1}{r^{1+2(h(t)-1/\alpha(t)})} \int_{\bbbr} \left( f(t+r,t,x) - f(t,t,x) \right)^2  m(dx) & = & \frac{1}{r} \int_{\bbbr} \mathbf{1}_{[t,t+r]}(x) dx\\
    & = & 1\\
    \end{eqnarray*}
thus (Cu14) holds.
\item (Cu15) $\forall v \in U$, $\forall u \in U$, 
	\begin{eqnarray*}
	 \frac{1}{|v-u|^2} \int_{\bbbr} \left| f(v,v,x) - f(v,u,x) \right|^2 m(dx) & = & 0 \\
	\end{eqnarray*}
	thus (Cu15) holds. \Box
 \end{itemize}

\noindent
{\bf Proof of theorem \ref{esplmmm}} 
 
\noindent
We want to apply theorem (\ref{vitesp}) with $f(t,u,x) = |t-x|^{H(u)-\frac{1}{\alpha(u)}} - |x|^{H(u)-\frac{1}{\alpha(u)}}$. Let us show that conditions (C1), (Cs2), (Cs3), (Cs4), (C5) and (C9) are satisfied.
 
\begin{itemize}
    \item (C1) The family of functions $u \to f(t,u,x)$ is differentiable for all $(u,t)$ in a neighbourhood of $t_0$ and almost all $x$ in $E$. The derivatives of $f$ with respect to $u$ read:
\begin{displaymath}
f'_u(t,w,x)=(h'(w)+\frac{\alpha'(w)}{\alpha^2(w)})\left[ (\log|t-x|)|t-x|^{h(w)-1/\alpha(w)} - (\log|x|)|x|^{h(w)-1/\alpha(w)}\right].
\end{displaymath}
    \item (Cs2) In \cite{FLV}, it is shown that, given $t_0 \in \bbbr$, one may choose $\varep>0$ small enough and numbers $a,b,h_{-},h_{+}$ with
$0<a <\alpha(w) < b<2 $, $0<h_{-} < h(w) < h_{+} <1$ and 
$\frac{a}{b}(\frac{1}{a}-\frac{1}{b})< h_{-}-(\frac{1}{a}-\frac{1}{b}) < h_{-}<h_{+}<h_{+}+ (\frac{1}{a}-\frac{1}{b}) <1-(\frac{1}{a}-\frac{1}{b})$
such that, for all $t$ and $w$ in $U:=(t_0-\varep, t_0+\varep)$ and all real $x$:
\begin{equation}\label{inegk}
|f(t,w,x)|,
|f'_{t_0}(t,w,x)| 
\leq  k_{1}(t,x)
\end{equation}

 where
\begin{equation}
k_{1}(t,x) = \left\{
\begin{array}{ll}
    c_{1}\max\{1, |t-x|^{h_{-}-1/a} +|x|^{h_{-}-1/a}\}
    & ( |x| \leq 1 + 2 \max_{t \in U}|t|) \\
    c_{2}|x|^{h_{+}-1/b-1} 
     & ( |x| > 1 + 2 \max_{t \in U}|t|) 
\end{array}
\right.
\end{equation}
for appropriately chosen constants $c_{1}$ and $c_{2}$.
One has, for all $\delta >0$, 
\begin{eqnarray*}
 \int_{\bbbr} \left[ \sup_{w \in U} |f(t,w,x)|^{\alpha(w)}\right]^{1+\delta} r(x)^{\delta} dx & \leq & \int_{\bbbr} \left( k_1(t,x)^a +k_1(t,x)^b \right)^{1+\delta} r(x)^{\delta} dx\\
 &\leq & K_{\delta} \int_{\bbbr} k_1(t,x)^{a(1+\delta)} r(x)^{\delta} dx \\
 & & + K_{\delta} \int_{\bbbr} k_1(t,x)^{b(1+\delta)} r(x)^{\delta} dx.\\
\end{eqnarray*}
Let us study $\int_{\bbbr} k_1(t,x)^{p(1+\delta)} r(x)^{\delta} dx$, where $p=a$ or $p=b$.
 
\begin{eqnarray*}
 \int_{\bbbr} k_1(t,x)^{p(1+\delta)} r(x)^{\delta} dx & = & \frac{\pi^{2\delta}}{3^{\delta}} \sum_{j=0}^{+\infty} (j+1)^{2 \delta} \int_{j}^{j+1} ( k_1(t,-x)^{p(1+\delta)} +k_1(t,x)^{p(1+\delta)})dx\\
 & = & \frac{\pi^{2\delta}}{3^{\delta}} \sum_{j=0}^{+\infty} (j+1)^{2 \delta} \int_{j}^{j+1} ( k_1(-t,x)^{p(1+\delta)} +k_1(t,x)^{p(1+\delta)})dx.\\ 
\end{eqnarray*}
We consider now $\int_{j}^{j+1} k_1(\pm t,x)^{p(1+\delta)}dx.$ There exists $K_{p,\delta} > 0$ such that, for all real $x$ such that $ |x| \leq 1 + 2 \max_{t \in U}|t|$: 
$$k_1(\pm t,x)^{p(1+\delta)} \leq K_{p,\delta} \left( 1 + |\pm t -x|^{p(1+\delta)(h_{-} -1/a)} + |x|^{p(1+\delta)(h_{-} -1/a)}\right),$$
and for all real $x$ such that $ |x| > 1 + 2 \max_{t \in U}|t|$, 
$$k_1(\pm t,x)^{p(1+\delta)} \leq K_{p,\delta} |x|^{p(1+\delta)(h_{+}-1/b-1)}.$$ 
Let $j_0 = \left[ 1 + 2 \max_{t \in U}|t|)\right].$
For $j <j_0$,
\begin{eqnarray*}
 \int_{j}^{j+1} k_1(\pm t,x)^{p(1+\delta)}  dx & \leq   & K_{p,\delta} (1 + \int_{j}^{j+1} |\pm t -x|^{p(1+\delta)(h_{-} -1/a)} dx + \int_{j}^{j+1} |x|^{p(1+\delta)(h_{-} -1/a)} dx).\\
 \end{eqnarray*}
Choose $\delta$ such that $p(1+\delta)(h_{-} -1/a)>-1$ (we show below that such a $\delta$ exists). Then
 
 \begin{eqnarray*}
 \int_{j}^{j+1} |\pm t -x|^{p(1+\delta)(h_{-} -1/a)} dx &= &  \int_{\pm t - j -1}^{\pm t -j} |u|^{p(1+\delta)(h_{-} -1/a)} du \\
 & \leq & \frac{ |\pm t -j |^{1+p(1+\delta)(h_{-} -1/a)} + |\pm t -j-1|^{1+p(1+\delta)(h_{-} -1/a)}}{1+p(1+\delta)(h_{-} -1/a)} \\
 & \leq & K_U |t|^{1+p(1+\delta)(h_{-} -1/a)} |1+j|^{1+p(1+\delta)(h_{-} -1/a)}\\
 & \leq & K_U (1+j)^{1+p(1+\delta)(h_{-} -1/a)}\\
\end{eqnarray*}
where $K_U>0$ depends on $U$ and may have changed from line to line. We deduce:
\begin{displaymath}
 \int_{j}^{j+1} k_1(\pm t,x)^{p(1+\delta)}  dx  \leq K_U(1+j^{1+p(1+\delta)(h_{-} -1/a)}).
\end{displaymath}

For $j=j_0$,
\begin{eqnarray*}
 \int_{j_0}^{j_0+1} k_1(\pm t,x)^{p(1+\delta)}  dx & \leq & K_U |j_0|^{1+p(1+\delta)(h_{-} -1/a)} + K_U \int_{j_0}^{j_0+1} |x|^{p(1+\delta)(h_{+}-1/b-1)} dx \\
 & \leq & K_U.
 \end{eqnarray*}
    
For $j>j_0$,
\begin{eqnarray*}
 \int_{j}^{j+1} k_1(\pm t,x)^{p(1+\delta)}  dx & \leq   & K_U  \int_{j}^{j+1} |x|^{p(1+\delta)(h_{+}-1/b-1)} dx\\
 & \leq & K_U j^{p(1+\delta)(h_{+}-1/b-1)}.\\
 \end{eqnarray*}
Finally,
\begin{eqnarray*}
 \sup_{t \in U} \int_{\bbbr} k_1(t,x)^{p(1+\delta)} r(x)^{\delta} dx & \leq & K_U \left(1+ \sum_{j=0}^{j_0-1} j^{2\delta}(1+j^{1+p(1+\delta)(h_{-} -1/a)} )  \right)\\
 & & + K_U \sum_{j=j_0+1}^{\infty} j^{2\delta + p(1+\delta)(h_{+}-1/b-1)}.\\
\end{eqnarray*}
To conclude, we need to show that we may chose $\delta > \frac{b}{a}-1$ such that $ p(1+\delta)(h_{-} -1/a)>-1$ and $2\delta + p(1+\delta)(h_{+}-1/b-1) < -1$, for $p=a$ and $p=b$. 
We consider several cases.\newline
\underline{First case} : $h_{-} -\frac{1}{a} \geq 0$ and $ h_{+}-\frac{1}{b} - 1 \leq -\frac{2}{a}$.

\smallskip

Let $\delta > \frac{b}{a}-1$. One has $p(1+\delta)(h_{-}-\frac{1}{a}) \geq 0 > -1$. We consider $1+2\delta + p(1+\delta)(h_{+}-1/b-1).$
 \begin{eqnarray*}
  1+2\delta + p(1+\delta)(h_{+}-1/b-1) & \leq & 1+2\delta-\frac{2}{a}p(1+\delta)\\
  & = & 1-\frac{2p}{a} + 2 \delta (1-\frac{p}{a}).
 \end{eqnarray*}
Since $1-\frac{2p}{a}<0$ and $1-\frac{p}{a} \leq 0$, $1+2\delta + p(1+\delta)(h_{+}-1/b-1) <0$.
\newline
\underline{Second case} : $h_{-} -\frac{1}{a} \geq 0$ and $ h_{+}-\frac{1}{b} - 1 > -\frac{2}{a}$.

\smallskip

Let $\delta \in \left(\frac{b}{a}-1 ,  \frac{\frac{1}{b}-\frac{1}{a} + 1 - h_{+}}{\frac{2}{a} -\frac{1}{b} -1 +h_{+}}\right)$. One has $p(1+\delta)(h_{-}-\frac{1}{a}) \geq 0 > -1$. We consider $1+2\delta + p(1+\delta)(h_{+}-1/b-1).$
\newline
\underline{For $p=a$} :
\begin{eqnarray*}
 1+2\delta + p(1+\delta)(h_{+}-1/b-1) & = & a \delta (\frac{2}{a} + h_{+}-\frac{1}{b}-1) + a (h_{+}-\frac{1}{b}-1+\frac{1}{a} )\\
 & < & a ( \frac{1}{b}-\frac{1}{a} + 1 - h_{+}) + a (h_{+}-\frac{1}{b}-1+\frac{1}{a} )\\
 & = & 0. 
 \end{eqnarray*}
 \underline{For $p=b$} :
 \begin{eqnarray*}
  1+2\delta + p(1+\delta)(h_{+}-1/b-1)& = & b \delta (\frac{1}{b} + h_{+}-1) + b (h_{+}-1).\\
 \end{eqnarray*}
If $\frac{1}{b} + h_{+}-1 \leq 0$, then $b \delta (\frac{1}{b} + h_{+}-1) + b (h_{+}-1) <0$. Else
\begin{eqnarray*}
 b \delta (\frac{1}{b} + h_{+}-1) + b (h_{+}-1) & < & b \frac{\frac{1}{b}-\frac{1}{a} + 1 - h_{+}}{\frac{2}{a} -\frac{1}{b} -1 +h_{+}} (\frac{1}{b} + h_{+}-1) + b (h_{+}-1)\\
 & = & \frac{b}{\frac{2}{a} -\frac{1}{b} -1 +h_{+}} (\frac{1}{a}-\frac{1}{b}) ( h_{+}-1 - \frac{1}{b})\\
 & < & 0.\\
\end{eqnarray*}

\underline{Third case} : $h_{-} -\frac{1}{a} < 0$ and $ h_{+}-\frac{1}{b} - 1 \leq -\frac{2}{a}$.
\smallskip

Let $\delta \in \left(\frac{b}{a}-1 ,  \frac{a h_{-}+ \frac{a}{b}-1}{1 -ah_{-}}\right)$.
\newline
\underline{For $p=a$} : 
\begin{eqnarray*}
 1+p(1+\delta)(h_{-}-\frac{1}{a}) & = & a h_{-}+\delta(ah_{-}-1)\\
 & > & a h_{-} + (ah_{-}-1)\frac{a h_{-}+ \frac{a}{b}-1}{1 -ah_{-}}\\
 & = & a h_{-} + 1 - \frac{a}{b} - a h_{-}\\
 & > & 0,\\
\end{eqnarray*}
and
\begin{eqnarray*}
 1+2\delta + p(1+\delta)(h_{+}-1/b-1) & = & a \delta (\frac{2}{a} + h_{+}-\frac{1}{b}-1) + a (h_{+}-\frac{1}{b}-1+\frac{1}{a} )\\
 & \leq & a (h_{+}-\frac{1}{b}-1+\frac{1}{a} )\\
 & \leq & -1\\
 & < & 0.\\
\end{eqnarray*}
\underline{For $p=b$} : 
\begin{eqnarray*}
 1+p(1+\delta)(h_{-}-\frac{1}{a}) & = & b ( h_{-}-\frac{1}{a} + \frac{1}{b}) + b \delta (h_{-}-\frac{1}{a} )\\
 & > & b ( h_{-}-\frac{1}{a} + \frac{1}{b}) + b (h_{-}-\frac{1}{a} ) \frac{a h_{-}+ \frac{a}{b}-1}{1 -ah_{-}}\\
 & = & b ( h_{-}-\frac{1}{a} + \frac{1}{b}) + b (\frac{1}{a} - \frac{1}{b} - h_{-})\\
 & = & 0,\\
\end{eqnarray*}
and
\begin{eqnarray*}
 1+2\delta + p(1+\delta)(h_{+}-1/b-1) & = & b \delta (\frac{1}{b} + h_{+}-1) + b (h_{+}-1)\\
 & \leq & b \delta (\frac{2}{b} - \frac{2}{a}) + b (h_{+}-1)\\
 & < & 0.\\
\end{eqnarray*}

\underline{Fourth case} : $h_{-} -\frac{1}{a} < 0$ and $ h_{+}-\frac{1}{b} - 1 > -\frac{2}{a}$.
\smallskip

Let $\delta \in \left(\frac{b}{a}-1 ,  \min(\frac{a h_{-}+ \frac{a}{b}-1}{1 -ah_{-}},\frac{\frac{1}{b}-\frac{1}{a} + 1 - h_{+}}{\frac{2}{a} -\frac{1}{b} -1 +h_{+}} ) \right)$.
\newline
\underline{For $p=a$} : 
\begin{eqnarray*}
 1+p(1+\delta)(h_{-}-\frac{1}{a}) & = & a h_{-}+\delta(ah_{-}-1)\\
 & > & a h_{-} + (ah_{-}-1)\frac{a h_{-}+ \frac{a}{b}-1}{1 -ah_{-}}\\
 & = & a h_{-} + 1 - \frac{a}{b} - a h_{-}\\
 & > & 0,\\
\end{eqnarray*}
and
\begin{eqnarray*}
 1+2\delta + p(1+\delta)(h_{+}-1/b-1) & = & a \delta (\frac{2}{a} + h_{+}-\frac{1}{b}-1) + a (h_{+}-\frac{1}{b}-1+\frac{1}{a} )\\
 & > & a (\frac{1}{b}-\frac{1}{a} + 1 - h_{+}) + a (h_{+}-\frac{1}{b}-1+\frac{1}{a} )\\
 & = & 0.\\
\end{eqnarray*}
\underline{For $p=b$} : 
\begin{eqnarray*}
 1+p(1+\delta)(h_{-}-\frac{1}{a}) & = & b ( h_{-}-\frac{1}{a} + \frac{1}{b}) + b \delta (h_{-}-\frac{1}{a} )\\
 & > & b ( h_{-}-\frac{1}{a} + \frac{1}{b}) + b (h_{-}-\frac{1}{a} ) \frac{a h_{-}+ \frac{a}{b}-1}{1 -ah_{-}}\\
 & = & b ( h_{-}-\frac{1}{a} + \frac{1}{b}) + b (\frac{1}{a} - \frac{1}{b} - h_{-})\\
 & = & 0,\\
\end{eqnarray*}
and
\begin{eqnarray*}
 1+2\delta + p(1+\delta)(h_{+}-1/b-1) & = & b \delta (\frac{1}{b} + h_{+}-1) + b (h_{+}-1).\\
\end{eqnarray*}
If $\frac{1}{b} + h_{+}-1 \leq 0$, then $ 1+2\delta + p(1+\delta)(h_{+}-1/b-1)< 0$, else
\begin{eqnarray*}
 b \delta (\frac{1}{b} + h_{+}-1) + b (h_{+}-1) & < & b (\frac{\frac{1}{b}-\frac{1}{a} + 1 - h_{+}}{\frac{2}{a} -\frac{1}{b} -1 +h_{+}})(\frac{1}{b} + h_{+}-1) + b (h_{+}-1)\\
 & = & \frac{b}{\frac{2}{a} -\frac{1}{b} -1 +h_{+}} (\frac{1}{a}-\frac{1}{b})(h_{+}-1-\frac{1}{b})\\
 & < & 0.\\
\end{eqnarray*}

    \item (Cs3) is obtained with (\ref{inegk}) for the same reason as in (Cs2).

    \item (Cs4) For $j$ large enough ($j > j^*$),
    \begin{eqnarray*}
 |f(t,w,x)\log(r(x))|^{\alpha(w)}  & \leq & K_1|k_1(t,x)|^{\alpha(w)}\\
& & + K_2 \sum_{j=j^*}^{+\infty} |f(t,w,x)|^{\alpha(w)} (\log(j))^d \mathbf{1}_{[-j,-j+1[ \cup [j-1,j[}(x).\\
\end{eqnarray*}
\begin{eqnarray*}
 |f(t,w,x)|^{\alpha(w)}\mathbf{1}_{[-j,-j+1[ \cup [j-1,j[}(x) &\leq& K_2 \frac{1}{|x|^{a(1
 +1/b-h_+)}} \mathbf{1}_{[-j,-j+1[ \cup [j-1,j[}(x)
\end{eqnarray*}
 ($K_2$ may have changed from line to line). Thus

\begin{eqnarray*}
\left[ \sup_{w \in U} \left[  \left|f(t,w,x) \log(r(x))  \right|^{\alpha(w)}  \right]\right]^{1+\delta} r(x)^{\delta} & \leq & K|k_1(t,x)|^{a(1+\delta)}r(x)^{\delta} + K|k_1(t,x)|^{b(1+\delta)}r(x)^{\delta}\\
& & + K \sum_{j=j^*}^{+\infty} \frac{j^{2\delta}(\log(j))^d}{|x|^{a(1+\delta)(1
 +1/b-h_+)}} \mathbf{1}_{[-j,-j+1[ \cup [j-1,j[}(x).
\end{eqnarray*}
Let $\delta > \frac{b}{a} - 1$ be such that (Cs2) holds. Since $2\delta +a(1+\delta)(h_{+}-1-\frac{1}{b}) <-1$, (Cs4) holds .

\item (C5) $X(t,u)$ (as a process in $t$) is localisable at $u$ with exponent $H(u) \in (h_-, h_+) \subset (0,1)$, with local form $X_u(t,u)$, and  $u \mapsto H(u)$ is a $C^1$ function (see \cite{LGLV}).

\item (C9)
\begin{eqnarray*}
 \frac{1}{r^{H(t)\alpha(t)}} \int_{\bbbr} | f(t+r,t,x)-f(t,t,x)|^{\alpha(t)} m(dx) & = &  \int_{\bbbr} \left| | 1-x |^{H(t)-\frac{1}{\alpha(t)}} - |x|^{H(t)-\frac{1}{\alpha(t)}}\right|^{\alpha(t)} dx \\
\end{eqnarray*}
so (C9) holds.

\end{itemize}

From theorem \ref{vitesp}, we obtain that 
 \begin{displaymath}
 \E \left[ |Y(t+\varep)-Y(t)|^{\eta}\right] \sim \varep^{\eta H(t)} \E \left[ |Y'_t(1)|^{\eta}\right].
\end{displaymath}
Since $Y'_t(1)$ is an $S_{\alpha(t)}(\sigma,0,0)$ random variable with $\sigma = \left( \int_{\bbbr} \left| | 1-x |^{H(t)-\frac{1}{\alpha(t)}} - |x|^{H(t)-\frac{1}{\alpha(t)}}\right|^{\alpha(t)} dx \right)^{\frac{1}{\alpha(t)}}$, property 1.2.17 of \cite{ST} allows to conclude. \Box

\noindent
{\bf Proof of theorem \ref{expolmmm}} 
 
\noindent
We want to apply theorems \ref{upbound} with $f(t,u,x) = |t-x|^{H(u)-\frac{1}{\alpha(u)}} - |x|^{H(u)-\frac{1}{\alpha(u)}}$ in order to obtain the inequality. Let us show that conditions (C6), (C7), (Cu8), (Cu10), (Cu11), (Cu12), (C13), (Cu14) and (Cu15) are satisfied.
\begin{itemize}
    \item (C6) Since $H(t) - \frac{1}{\alpha(t)} \geq 0$, (C6) holds.
    \item (C7) We also use the fact that $H(t) - \frac{1}{\alpha(t)} \geq 0$ in order to prove that (C7) holds.
    \item (Cu8) $\forall v \in U$, $\forall u \in U$, $\forall x \in \bbbr$,
    \begin{eqnarray*}
    \frac{1}{|v-u|^{h(u)-1/ \alpha(u)}}\left| f(v,u,x) -f(u,u,x) \right|  & = & \frac{1}{|v-u|^{H(u)-1/ \alpha(u)}} \left| |v-x|^{H(u)-\frac{1}{\alpha(u)}} - |u-x|^{H(u)-\frac{1}{\alpha(u)}} \right|\\
    & \leq & 1\\
    \end{eqnarray*}
    thus (Cu8) holds.
    \item (Cu10) $\forall v \in U$,$\forall u \in U$,
    \begin{eqnarray*}
     \frac{1}{|v-u|^{1+p(h(u)-\frac{1}{\alpha(u)})}} \int_{\bbbr} \left| f(v,u,x) - f(u,u,x) \right|^p m(dx) & = & \int_{\bbbr} \left| | 1-x |^{H(u)-\frac{1}{\alpha(u)}} - |x|^{H(u)-\frac{1}{\alpha(u)}}\right|^{p} dx\\
 \end{eqnarray*}
 so (Cu10) holds.
    \item (Cu11) $\forall  v \in U$, $\forall  u \in U$, 
	\begin{displaymath}
	\int_{\bbbr} \left| f(v,u,x) \right|^2 m(dx) = v^{1+2(H(u)-\frac{1}{\alpha(u)})} \int_{\bbbr} \left| | 1-x |^{H(u)-\frac{1}{\alpha(u)}} - |x|^{H(u)-\frac{1}{\alpha(u)}}\right|^{2} dx
	\end{displaymath}
	thus (Cu11) holds.
    \item (Cu12) For the same reason as (Cu11), (Cu12) holds.
    \item (C13) For $t \neq 0$, one can choose $U$ such that $\inf_{v \in U} v^{1+2(H(v)-\frac{1}{\alpha(v)})} > 0$ thus (C13) holds.
    \item (Cu14) 
    \begin{eqnarray*}
    \frac{1}{r^{1+2(h(t)-1/\alpha(t)})} \int_{\bbbr} \left( f(t+r,t,x) - f(t,t,x) \right)^2  m(dx) & = & \int_{\bbbr} \left| | 1-x |^{H(t)-\frac{1}{\alpha(t)}} - |x|^{H(t)-\frac{1}{\alpha(t)}}\right|^{2} dx\\
    \end{eqnarray*}
thus, choosing $g(t)= \int_{\bbbr} \left| | 1-x |^{H(t)-\frac{1}{\alpha(t)}} - |x|^{H(t)-\frac{1}{\alpha(t)}}\right|^{2} dx$, (Cu14) holds.
\item (Cu15) $\forall v \in U$, $\forall u \in U$, 
	\begin{eqnarray*}
	 \frac{1}{|v-u|^2} \int_{\bbbr} \left| f(v,v,x) - f(v,u,x) \right|^2 m(dx) & = &\\ 
	  \frac{1}{|v-u|^2} \int_{\bbbr} \left| | v-x |^{H(v)-\frac{1}{\alpha(v)}}  - | v-x |^{H(u)-\frac{1}{\alpha(u)}} - | x |^{H(v)-\frac{1}{\alpha(v)}} + | x |^{H(u)-\frac{1}{\alpha(u)}}\right|^2 dx & &\\
	\end{eqnarray*}
	thus (Cu15) holds \Box
 \end{itemize}
 
\section{Proof of Theorem \ref{expolevyexact}}\label{proolevy}

Recall the definition of the L\'evy Multistable field on $[0,1]$:
\begin{displaymath}
X(v,u)= C^{1/\alpha(u)}_{\alpha(u)} \sum_{i=1}^{\infty} \gamma_i \Gamma_i^{-1/\alpha(u)} \mathbf{1}_{[0,v]}(V_i).
\end{displaymath}

To prove Theorem \ref{expolevyexact}, we need a series of lemma:
 
\begin{lem}\label{leminegtayl} Assume $\alpha$ is $\mathcal{C}^1$. Then, for all $u \in (0,1)$, almost surely,
\begin{displaymath}
\sup_{v \in [0,1]}\limits  \frac{|X(v,v)-X(v,u) |}{|v-u|}< + \infty.
\end{displaymath}   
\end{lem}

\noindent
{\bf Proof} 

\noindent
in the case of the Lévy multistable field, \eqref{decompchamp} reads:
\begin{displaymath}
X(v,v)-X(v,u) = (v-u) \left(  \sum_{i=1}^{+\infty} Z_i^1(v) +\sum_{i=1}^{+\infty} Z_i^3(v) +\sum_{i=1}^{+\infty} Y_i^1(v) +\sum_{i=1}^{+\infty} Y_i^3(v)\right),
\end{displaymath}
where $Z_i^1,\ldots$ are defined as above.
Let $A >0$ and $B>0$ be constants such that $\forall w \in U$, $|a'(w) | \leq A$ and $ | a(w)\frac{\alpha'(w)}{\alpha^2(w)}|\leq B$.
We write $\sum_{i=1}^{+\infty}\limits Z_i^1(v)= \sum_{j=1}^{+\infty}\limits \left( \sum_{i=2^j}^{2^{j+1}-1}\limits Z_i^1(v)\right) =: \sum_{j=1}^{+\infty}\limits X_j^1(v) $ and $\sum_{i=1}^{+\infty}\limits Z_i^3(v) = \sum_{j=1}^{+\infty}\limits \left( \sum_{i=2^j}^{2^{j+1}-1}\limits Z_i^3(v)\right) =: \sum_{j=1}^{+\infty}\limits X_j^3(v) .$ We consider $\liminf_j \{ \sup_{v \in [0,1]}\limits | X_j^1(v)| \leq \frac{ A j\sqrt{2^j} }{2^{j/d}} \}$ and $\liminf_j \{ \sup_{v \in [0,1]}\limits | X_j^3(v)| \leq \frac{\log(2) B j(j+1) \sqrt{2^j}}{2^{j/d}} \}$. Let $V^{(1)}, V^{(2)}, ... ,V^{(2^j)}$ denote the order statistics of the
$V_i$ ({\it i.e.}$V^{(1)}=\min V_i$, \ldots). Then:
\begin{displaymath}
\{\sup_{v \in [0,1]}\limits | X_j^1(v)| > \frac{ A j\sqrt{2^j} }{2^{j/d}} \} \subset \cup_{ N \geq 1}^{2^j} \cup_{l_1,..., l_N \in \llbracket 2^j,2^{j+1}-1\rrbracket} \left( \{ |\sum_{i = 1}^{N} \gamma_{l_i} a'(w_{l_i}) l_i^{-1/ \alpha(w_{l_i})} |  > \frac{ A j\sqrt{2^j} }{2^{j/d}} \} \right...
\end{displaymath}
\begin{displaymath}
 \left... \cap \{ V^{(1)} = V_{l_1}, V^{(2)} = V_{l_2} , ... , V^{(N)} = V_{l_N}\} \right).
\end{displaymath}
\begin{eqnarray*}
\P \left(\sup_{v \in [0,1]}\limits | X_j^1(v)| > \frac{ A j\sqrt{2^j} }{2^{j/d}}  \right) & \leq & \sum_{N=1}^{2^{j}} \sum_{l_1,..., l_N \in \llbracket 2^j,2^{j+1}-1\rrbracket} \frac{(2^j -N)!}{(2^j)!}\P \left(  | \sum_{i = 1}^{N} \gamma_{l_i} a'(w_{l_i}) l_i^{-1/ \alpha(w_{l_i})}| >\frac{ A j\sqrt{2^j} }{2^{j/d}} \right) \\
& \leq & \sum_{N=1}^{2^{j}} \sum_{l_1,..., l_N \in \llbracket 2^j,2^{j+1}-1\rrbracket}\frac{(2^j -N)!}{(2^j)!}\P \left(  | \sum_{i = 1}^{N} \gamma_{l_i} \frac{a'(w_{l_i})}{A} \frac{2^{j/d}}{l_i^{1/ \alpha(w_{l_i})}}| > j\sqrt{N}\right) \\
& \leq & \sum_{N=1}^{2^{j}} \sum_{l_1,..., l_N \in \llbracket 2^j,2^{j+1}-1\rrbracket}\frac{(2^j -N)!}{(2^j)!} 2 e^{-\frac{j^2}{2}}\\
& \leq & 2 e^{-\frac{j^2}{2}} \sum_{N=1}^{2^{j}} \frac{1}{N!}\\
& \leq & 2e^{1-\frac{j^2}{2}}\\
\end{eqnarray*}
where we have used the following inequality (lemma 1.5, chapter 1 in \cite{LTal}):
\begin{displaymath}
 \P \left( |\sum_{i=1}^n u_i| \geq \lambda \sqrt{n}\right) \leq 2 e^{-\frac{\lambda^2}{2}}
\end{displaymath}
for $(u_i)_i$ independent centered random variables verifying
$-1 \leq u_i \leq 1$, with $u_i = \gamma_{l_i} \frac{a'(w_{l_i})}{A} \frac{2^{j/d}}{l_i^{1/ \alpha(w_{l_i})}} $ 
and $\lambda=j$.

We deduce that $\P \left(\liminf_j \left\{ \sup_{v \in [0,1]}\limits | X_j^1(v)| \leq \frac{A j \sqrt{2^j} }{2^{j/d}} \right\} \right) = 1$. 

Similarly:
\begin{eqnarray*}
\P \left(\sup_{v \in [0,1]}\limits | X_j^3(v)| > \frac{\log(2) B j(j+1) \sqrt{2^j}}{2^{j/d}} \right)& \leq & 2e^{1-\frac{j^2}{2}}\\
\end{eqnarray*}
and $\P \left(\liminf_j \left\{ \sup_{v \in [0,1]}\limits | X_j^3(v)| \leq \frac{\log(2) B j(j+1) \sqrt{2^j}}{2^{j/d}} \right\}\right) = 1 $.
We work on the event 
\begin{displaymath}
\liminf_j \left\{ \sup_{v \in [0,1]}\limits | X_j^1(v)| \leq \frac{A j \sqrt{2^j} }{2^{j/d}} \right\} \cap  \liminf_j \left\{ \sup_{v \in [0,1]}\limits | X_j^3(v)| \leq \frac{\log(2) B j(j+1) \sqrt{2^j}}{2^{j/d}} \right\} \cap \liminf_i \left\{\Gamma_i >1 \right\}.
\end{displaymath}
There exists $J_0 \in \mathbb{N}$ such that $\forall j \geq J_0$, $\sup_{v \in [0,1]}\limits | X_j^1(v)| \leq  \frac{A j \sqrt{2^j} }{2^{j/d}} $ and $\sup_{v \in [0,1]}\limits | X_j^3(v)| \leq  \frac{\log(2) B j(j+1) \sqrt{2^j}}{2^{j/d}} $.
\begin{displaymath}
\left| \sum_{i=1}^{+\infty} Z_i^1(v)\right | \leq \sum_{j=0}^{2^{J_0}-1}\frac{A}{i^{1/d}} + \sum_{j=J_0}^{+\infty} A \frac{j}{2^{j(\frac{1}{d}-\frac{1}{2})}}
\end{displaymath}
and
\begin{displaymath}
\left| \sum_{i=1}^{+\infty} Z_i^3(v)\right | \leq \sum_{j=0}^{2^{J_0}-1}\frac{B \log(i)}{i^{1/d}} + \sum_{j=J_0}^{+\infty} B\log(2) \frac{j(j+1)}{2^{j(\frac{1}{d}-\frac{1}{2})}},
\end{displaymath}
thus $\sup_{v \in [0,1]}\limits \left| \sum_{i=1}^{+\infty}\limits Z_i^1(v)\right| < +\infty $ and $\sup_{v \in [0,1]}\limits \left|\sum_{i=1}^{+\infty}\limits Z_i^3(v) \right| < + \infty$.\par
Fix $i_0 \in \mathbb{N}$ such that $\forall i \geq i_0$, $\Gamma_i > 1$.  
\begin{displaymath}
|\sum_{i=1}^{i_0} Y_i^1(v)| \leq A \sum_{i=1}^{i_0}( \frac{1}{\Gamma_i^{1/c}} + \frac{1}{i^{1/d}})
\end{displaymath}
and
\begin{displaymath}
|\sum_{i=1}^{i_0} Y_i^3(v)| \leq B \sum_{i=1}^{i_0} \left( | \frac{\log \Gamma_i}{\Gamma_i^{1/c}}| + \frac{\log(i)}{i^{1/d}}\right).
\end{displaymath}
\begin{eqnarray*}
|\sum_{i=i_0}^{+\infty} Y_i^1(v)| & \leq & A \sum_{i=i_0}^{+\infty} | \Gamma_i^{-1/\alpha(x_i)} - i^{-1/\alpha(x_i)}|\mathbf{1}_{ \{1< \Gamma_i \leq \frac{i}{2} \} }\\
&  &+ A \sum_{i=i_0}^{+\infty} | \Gamma_i^{-1/\alpha(x_i)} - i^{-1/\alpha(x_i)}|\mathbf{1}_{ \{\frac{i}{2}< \Gamma_i \leq 2i \} } \\
& & +A \sum_{i=i_0}^{+\infty} | \Gamma_i^{-1/\alpha(x_i)} - i^{-1/\alpha(x_i)}|\mathbf{1}_{ \{ \Gamma_i > 2i \} },\\
\end{eqnarray*}

\begin{eqnarray*}
|\sum_{i=i_0}^{+\infty} Y_i^1(v)| & \leq & 2A \sum_{i=i_0}^{+\infty} (\mathbf{1}_{ \{1< \Gamma_i \leq \frac{i}{2} \} } +\mathbf{1}_{ \{ \Gamma_i > 2i \} } )+ A \sum_{i=i_0}^{+\infty} | \Gamma_i^{-1/\alpha(x_i)} - i^{-1/\alpha(x_i)}|\mathbf{1}_{ \{\frac{i}{2}< \Gamma_i \leq 2i \} } \\
& \leq & 2A \sum_{i=i_0}^{+\infty} (\mathbf{1}_{ \{1< \Gamma_i \leq \frac{i}{2} \} } +\mathbf{1}_{ \{ \Gamma_i > 2i \} } ) + K_{c,d} \sum_{i=i_0}^{+\infty} \frac{1}{i^{\frac{1}{d}}} | \frac{\Gamma_i}{i}-1|.
\end{eqnarray*}

\begin{eqnarray*}
|\sum_{i=i_0}^{+\infty} Y_i^3(v)| & \leq &  B \sum_{i=i_0}^{+\infty} | \log(\Gamma_i) \Gamma_i^{-1/\alpha(x_i)} - \log(i)i^{-1/\alpha(x_i)}|\mathbf{1}_{ \{1< \Gamma_i \leq \frac{i}{2} \} } \\
& &+ B \sum_{i=i_0}^{+\infty} | \log(\Gamma_i)\Gamma_i^{-1/\alpha(x_i)} - \log(i) i^{-1/\alpha(x_i)}|\mathbf{1}_{ \{\frac{i}{2}< \Gamma_i \leq 2i \} }\\
& &  +B \sum_{i=i_0}^{+\infty} | \log(\Gamma_i) \Gamma_i^{-1/\alpha(x_i)} - \log(i) i^{-1/\alpha(x_i)}|\mathbf{1}_{  \{\Gamma_i > 2i \} },\\
\end{eqnarray*}

\begin{eqnarray*}
|\sum_{i=i_0}^{+\infty} Y_i^3(v)| & \leq & K \sum_{i=i_0}^{+\infty} \log(i) (\mathbf{1}_{ \{1< \Gamma_i \leq \frac{i}{2} \} } +\mathbf{1}_{  \{\Gamma_i > 2i \} } )\\
& &+ B \sum_{i=i_0}^{+\infty} | \log(\Gamma_i) \Gamma_i^{-1/\alpha(x_i)} - \log(i) i^{-1/\alpha(x_i)}|\mathbf{1}_{ \{\frac{i}{2}< \Gamma_i \leq 2i \} } \\
& \leq & K \sum_{i=i_0}^{+\infty} \log(i) (\mathbf{1}_{ \{1< \Gamma_i \leq \frac{i}{2} \} } +\mathbf{1}_{  \{\Gamma_i > 2i \} } ) + K_{c,d} \sum_{i=i_0}^{+\infty} \frac{\log(i)}{i^{\frac{1}{d}}} | \frac{\Gamma_i}{i}-1|.
\end{eqnarray*}
Finally,  $\sup_{v \in [0,1]}\limits \left| \sum_{i=1}^{+\infty}\limits Y_i^1(v)\right| < +\infty $ and $\sup_{v \in [0,1]}\limits \left|\sum_{i=1}^{+\infty}\limits Y_i^3(v) \right| < + \infty$.\par
As a consequence, $\sup_{v \in [0,1]}\limits  \frac{|X(v,v)-X(v,u) |}{|v-u|}< + \infty$ \Box

\begin{lem}\label{lemexpostab}
For all $u \in (0,1)$ and all $\eta \in (0, \frac{1}{\alpha(u)})$, one has, almost surely,
\begin{displaymath}
\sup_{v \in [0,1]} \left| \frac{X(v,u)-X(u,u)}{|v-u|^{\eta}}\right| < +\infty.
\end{displaymath}
\end{lem}

\noindent
{\bf Proof} 

\noindent
Let $\eta \in (0, \frac{1}{\alpha(u)})$, $m \in \mathbb{N}$ , $C_j=\cap_{i=2^j}^{2^{j+1}-1} \{ V_i \not\in [u-\frac{1}{j^2 2^j},u+\frac{1}{j^2 2^j} ]\}$,
\begin{displaymath}
D_j^m=\left\{ \sup_{\frac{1}{2^{m+1}} \leq  |v -u| \leq \frac{1}{2^m} }\left|\sum_{i=2^j}^{2^{j+1}-1}\limits \gamma_i i^{-1/ \alpha(u)} \frac{\mathbf{1}_{[u,v]}(V_i)}{|v-u |^{\eta}} \right| \leq \frac{1}{j^2} \right\},
\end{displaymath}
and $D_j = \cap_{m\geq 0} D_j^m $. $D_j$ may be written:

\begin{displaymath}
D_j=\left\{ \sup_{v \in [0,1]} \left|\sum_{i=2^j}^{2^{j+1}-1}\limits \gamma_i i^{-1/ \alpha(u)} \frac{\mathbf{1}_{[u,v]}(V_i)}{|v-u |^{\eta}} \right| \leq \frac{1}{j^2} \right\}.
\end{displaymath}
Let us evaluate $\liminf C_j$.
\begin{displaymath}
\P \left( \overline{C_j}\right) \leq \sum_{i=2^{j}}^{2^{j+1}-1} \frac{1}{j^2 2^j} = \frac{1}{j^2}
\end{displaymath}
and thus $\P \left( \liminf_j C_j\right) = 1$. Now,
\begin{eqnarray*}
\P \left( \overline{D_j}\right) &\leq  & \frac{1}{j^2} + \P ( \overline{D_j}\cap C_j)\\
& = & \frac{1}{j^2}+ \P \left( \cup_{m \geq 0} (\overline{D_j^m}\cap C_j) \right)\\
& \leq & \frac{1}{j^2}+ \sum_{m =0}^{+ \infty} \P \left( \overline{D_j^m}\cap C_j\right).\\
\end{eqnarray*}
\newline
We consider several cases, depending on the respective values of $j$ and $m$:
\begin{itemize}
\item If $m > j+ \frac{2}{\log (2)} \log j$,
\end{itemize}
\begin{displaymath}
\P \left( \overline{D_j^m}\cap C_j\right) = 0.
\end{displaymath}
\begin{itemize}
\item If $ j+ \frac{2}{\log (2)} \log j \geq m \geq j$,
\end{itemize}

\begin{eqnarray*}
\P \left( \overline{D_j^m}\right)& \leq &  \P \left(\sup_{\frac{1}{2^{m+1}} \leq  |v -u| \leq \frac{1}{2^m} }\left|\sum_{i=2^j}^{2^{j+1}-1}\limits \gamma_i i^{-1/ \alpha(u)} \mathbf{1}_{[u,v]}(V_i) \right| \geq \frac{1}{2^{(m+1)\eta}j^2}\right).\\
\end{eqnarray*}
Let $J_0 \in \mathbb{N}$ be such that for all $ j > J_0$, $ 2^{j(\frac{1}{\alpha(u)}- \eta)}>2^{\eta}j^{3+\frac{2\eta}{\log(2)}}$. The event: 
\begin{displaymath}
\{\sup_{\frac{1}{2^{m+1}} \leq  |v -u| \leq \frac{1}{2^m} }\left|\sum_{i=2^j}^{2^{j+1}-1}\limits \gamma_i i^{-1/ \alpha(u)} \mathbf{1}_{[u,v]}(V_i) \right| \geq \frac{1}{2^{(m+1)\eta}j^2} \}
\end{displaymath}
is included in the event
\begin{displaymath}
\cup_{ N \geq 1}^{2^j} \left( \cup_{l_1,..., l_N \in \llbracket 2^j,2^{j+1}-1\rrbracket}  \{ |\sum_{i = 1}^{N} \gamma_{l_i}  l_i^{-1/ \alpha(u)} |  >  \frac{1}{2^{(m+1)\eta}j^2} \} \cap  ( \cap_{i=1}^N \{ |V_{l_i} - u | \in [\frac{1}{2^{m+1}}, \frac{1}{2^{m}} ] \} ) \right...
\end{displaymath}
\begin{displaymath}
\left ...\cap  ( \cap_{k \neq l_i}\{ |V_{k} - u | \notin [\frac{1}{2^{m+1}}, \frac{1}{2^{m}} ] \}) \right).
\end{displaymath}
Notice that for $j \geq J_0$ and $N <j$, $\P \left(|\sum_{i = 1}^{N}\limits \gamma_{l_i}  l_i^{-1/ \alpha(u)} |  >  \frac{1}{2^{(m+1)\eta}j^2} \right) = 0$, and thus
\begin{eqnarray*}
\P \left( \overline{D_j^m}\right)& \leq & \sum_{N=j}^{2^j} \sum_{l_1,..., l_N \in \llbracket 2^j,2^{j+1}-1\rrbracket} \P \left(|\sum_{i = 1}^{N} \gamma_{l_i}  l_i^{-1/ \alpha(u)} |  >  \frac{1}{2^{(m+1)\eta}j^2} \right) \P \left( \cap_{i=1}^N \{ |V_{l_i} - u | \in [\frac{1}{2^{m+1}}, \frac{1}{2^{m}} ] \}\right)\\
& \leq & \sum_{N=j}^{2^j} \frac{1}{2^{(m+1)N}} \sum_{l_1,..., l_N \in \llbracket 2^j,2^{j+1}-1\rrbracket} \P \left(|\sum_{i = 1}^{N} \gamma_{l_i}  l_i^{-1/ \alpha(u)} |  >  \frac{1}{2^{(m+1)\eta}j^2} \right)\\
& \leq & \sum_{N=j}^{2^j} \frac{1}{2^{(m+1)N}} \sum_{l_1,..., l_N \in \llbracket 2^j,2^{j+1}-1\rrbracket} j^4 2^{2(m+1)\eta} \sum_{i=2^j}^{2^{j+1}-1} \frac{1}{i^{\frac{2}{\alpha(u)}}}\\
& \leq & \sum_{N=j}^{2^j} \frac{j^4 2^{2(m+1)\eta}}{2^{(m+1)N}} 2^{j(1-\frac{2}{\alpha(u)})} C_{2^j}^{N}\\
& \leq & j^4 2^{2(j+ \frac{2}{\log (2)} \log j+1)\eta -j\frac{2}{\alpha(u)}}\sum_{N=j}^{2^j} \frac{2^j C_{2^j}^{N}}{2^{(m+1)N}}\\
&  \leq & j^{4+ \frac{4 \eta}{\log(2)}} 2^{2j(\eta - \frac{1}{\alpha(u)})} \sum_{N=j}^{2^j} \frac{2^{j-N} 2^{(j-m)N}}{N!}\\
& \leq &  3 j^{4+ \frac{4 \eta}{\log(2)}} 2^{2j(\eta - \frac{1}{\alpha(u)})} .\\
\end{eqnarray*}
\begin{itemize}
\item When $ j \geq  m \geq \frac{\log(j)}{\log(2)}$, the same computations lead to: 
\end{itemize}
\begin{eqnarray*}
 & &\sum_{N=j2^{j-m}}^{2^j} \sum_{l_1,..., l_N \in \llbracket 2^j,2^{j+1}-1\rrbracket} \P \left(|\sum_{i = 1}^{N} \gamma_{l_i}  l_i^{-1/ \alpha(u)} |  >  \frac{1}{2^{(m+1)\eta}j^2} \right) \P \left( \cap_{i=1}^N \{ |V_{l_i} - u | \in [\frac{1}{2^{m+1}}, \frac{1}{2^{m}} ] \}\right)\\
 & \leq & \sum_{N=j2^{j-m}}^{2^j} \frac{j^4 2^{2(m+1)\eta}}{2^{(m+1)N}} 2^{j(1-\frac{2}{\alpha(u)})} C_{2^j}^{N}\\
 & \leq & j^4 2^{2(m+1)\eta -2j/\alpha(u)}\sum_{N=j2^{j-m}}^{2^j} \frac{2^{j-N} 2^{(j-m)N}}{N!}\\
 & \leq & j^4 2^{2\eta}2^{2j(\eta -\frac{1}{\alpha(u)})}\sum_{N=j2^{j-m}}^{+\infty} \frac{2^{(j-m)N}}{N!}\\
 & \leq & K j^4 2^{2j(\eta -\frac{1}{\alpha(u)})} \frac{e^{2^{j-m}} 2^{(j-m)(j2^{j-m}+1)}}{(j2^{j-m}+1)!}\\
\end{eqnarray*}
where we have used the estimate $\sum_{n \geq N} \frac{x^n}{n!} \leq e^x \frac{x^{N+1}}{(N+1)!}$.
We arrive at:
\begin{eqnarray*}
\P \left( \overline{D_j^m}\right)& \leq & K j^4 2^{2j(\eta -\frac{1}{\alpha(u)})} + \\
& &\sum_{N=1}^{j2^{j-m}} \frac{1}{2^{(m+1)N}} (1-\frac{1}{2^{m+1}})^{2^j-N}\sum_{l_1,..., l_N \in \llbracket 2^j,2^{j+1}-1\rrbracket} \P \left(|\sum_{i = 1}^{N} \gamma_{l_i}  l_i^{-1/ \alpha(u)} |  >  \frac{1}{2^{(m+1)\eta}j^2} \right). \\
\end{eqnarray*}
We need to distinguish two cases depending on the value of $\eta$. If $\eta \leq \frac{1}{2}$, fix $J_1 \in \mathbb{N}$ such that for all $j \geq J_1$, $ 2^{j(\frac{1}{\alpha(u)}-\frac{1}{2})} > 2^{1/ \alpha(u)}j^3 \sqrt{j}$. If $\eta > \frac{1}{2}$, fix $J_1 \in \mathbb{N}$ such that for all $j \geq J_1$, $2^{j(\frac{1}{\alpha(u)}-\eta)} > 2^{1/ \alpha(u)}j^3 \sqrt{j}$. Then for all $\eta$ and all $j \geq J_1$, one has $\frac{2^{j/ \alpha(u)}}{j^3 \sqrt{j2^{j-m}}2^{(m+1)\eta}} \geq 1$ and 
\begin{eqnarray*}
\P \left(|\sum_{i = 1}^{N} \gamma_{l_i}  l_i^{-1/ \alpha(u)} |  >  \frac{1}{2^{(m+1)\eta}j^2} \right) & \leq & \P \left(\left|\sum_{i = 1}^{N} \gamma_{l_i}  \left(\frac{2^j}{l_i}\right)^{1/ \alpha(u)} \right|  >  j \sqrt{N} \right)\\
& \leq & 2e^{-j^2/2}.\\
\end{eqnarray*}
We then get
\begin{eqnarray*}
\P \left( \overline{D_j^m}\right)& \leq & K j^4 2^{2j(\eta -\frac{1}{\alpha(u)})} + \sum_{N=1}^{j2^{j-m}} \frac{1}{2^{(m+1)N}} (1-\frac{1}{2^{m+1}})^{2^j-N} C_{2^j}^{N} 2e^{-j^2/2}\\
& \leq &  K j^4 2^{2j(\eta -\frac{1}{\alpha(u)})} + 2e^{-j^2/2} \sum_{N=1}^{2^{j}} \frac{1}{2^{(m+1)N}} (1-\frac{1}{2^{m+1}})^{2^j-N} C_{2^j}^{N}\\
& \leq &  K j^4 2^{2j(\eta -\frac{1}{\alpha(u)})} + 2e^{-j^2/2}.\\
\end{eqnarray*}
\begin{itemize}
\item Assume finally that $  m \leq \frac{\log(j)}{\log(2)}$.
\end{itemize}
Fix $J_2 \in \mathbb{N}$ such that for all $j \geq J_2$, $2^{j(\frac{1}{\alpha(u)}-\frac{1}{2})} > 2^{1/ \alpha(u)} j^{3+\eta}$. Then, for $j \geq J_2$, one has $\frac{2^{j/ \alpha(u)}}{j^3 \sqrt{2^j}2^{(m+1)\eta}} \geq 1$ and computations similar the ones above lead to 
\begin{eqnarray*}
\P \left( \overline{D_j^m}\right)& \leq & \sum_{N=1}^{2^{j}}\frac{1}{2^{(m+1)N}} (1-\frac{1}{2^{m+1}})^{2^j-N}\sum_{l_1,..., l_N \in \llbracket 2^j,2^{j+1}-1\rrbracket} \P \left(|\sum_{i = 1}^{N} \gamma_{l_i}  l_i^{-1/ \alpha(u)} |  >  \frac{1}{2^{(m+1)\eta}j^2} \right)\\
& \leq & \sum_{N=1}^{2^{j}}\frac{1}{2^{(m+1)N}} (1-\frac{1}{2^{m+1}})^{2^j-N}\sum_{l_1,..., l_N \in \llbracket 2^j,2^{j+1}-1\rrbracket} \P \left(|\sum_{i = 1}^{N} \gamma_{l_i}  (\frac{2^{j/ \alpha(u)}}{l_i})^{1/ \alpha(u)} |  >  j\sqrt{N} \right)\\
& \leq & 2e^{-j^2/2}.\\
\end{eqnarray*}
We thus get that, for $j \geq \max (J_0,J_1,J_2)$, 
\begin{eqnarray*}
\sum_{m =0}^{+ \infty} \P \left( \overline{D_j^m}\cap C_j\right) & \leq K \log(j)j^{4+ \frac{4 \eta}{\log(2)}} 2^{2j(\eta - \frac{1}{\alpha(u)})},
\end{eqnarray*}

and thus $\P \left( \liminf_j D_j\right) = 1$. 

On the event $\liminf_j C_j \cap \liminf_j D_j$, we may fix $j_0 \in \mathbb{N}$ such that for all $j \geq j_0$,
\begin{displaymath}
\sup_{v \in [0,1]}\limits \left|\sum_{i=2^j}^{2^{j+1}-1}\limits \gamma_i i^{-1/ \alpha(u)} \frac{\mathbf{1}_{[u,v]}(V_i)}{|v-u |^{\eta}} \right| \leq \frac{1}{j^2}.
\end{displaymath}
Since $\sup_{v \in [0,1]}\limits \left|\sum_{i=1}^{2^{j_0}-1}\limits \gamma_i i^{-1/ \alpha(u)} \frac{\mathbf{1}_{[u,v]}(V_i)}{|v-u |^{\eta}} \right|< +\infty$, we obtain
\begin{displaymath}
\sup_{v \in [0,1]}\limits \left|\sum_{i=1}^{+\infty}\limits \gamma_i i^{-1/ \alpha(u)} \frac{\mathbf{1}_{[u,v]}(V_i)}{|v-u |^{\eta}} \right|< +\infty.
\end{displaymath}

Let us now deal with
\begin{displaymath}
E_j=\left\{ \sup_{v \in [0,1]}\limits \left|\sum_{i=2^j}^{2^{j+1}-1}\limits \gamma_i (\Gamma_i^{-1/ \alpha(u)} - i^{-1/ \alpha(u)}) \frac{\mathbf{1}_{[u,v]}(V_i)}{|v-u |^{\eta}} \right| \leq \frac{1}{j^2}\right\}.
\end{displaymath}
\begin{eqnarray*}
\P \left( \overline{E_j}\right) &\leq  & \frac{1}{j^2} + \P ( \overline{E_j}\cap C_j)\\
& \leq & \frac{1}{j^2} + \P \left( 2^{j\eta}j^{2\eta} \sup_{v \in [0,1]}\limits \left|\sum_{i=2^j}^{2^{j+1}-1}\limits \gamma_i (\Gamma_i^{-1/ \alpha(u)} - i^{-1/ \alpha(u)}) \mathbf{1}_{[u,v]}(V_i) \right|> \frac{1}{j^2} \right)\\
& \leq & \frac{1}{j^2} + \P \left( \sum_{i=2^j}^{2^{j+1}-1}\limits \left| (\Gamma_i^{-1/ \alpha(u)} - i^{-1/ \alpha(u)}) \right|> \frac{1}{ 2^{j\eta}j^{2(1+ \eta)}} \right) \\
& \leq & \frac{1}{j^2} + 2^{j\eta}j^{2(1+ \eta)} \sum_{i=2^j}^{2^{j+1}-1}\limits \E | \Gamma_i^{-1/ \alpha(u)} - i^{-1/ \alpha(u)}|\\
& \leq & \frac{1}{j^2} + 2^{j\eta}j^{2(1+ \eta)} \sum_{i=2^j}^{2^{j+1}-1}\limits 2(\P (\Gamma_i < \frac{i}{2}) + \P(\Gamma_i >2i) ) \\
& &+ 2^{j\eta}j^{2(1+ \eta)} \sum_{i=2^j}^{2^{j+1}-1}\limits \E |\Gamma_i^{-1/ \alpha(u)} - i^{-1/ \alpha(u)}|\mathbf{1}_{\{ \frac{i}{2} < \Gamma_i < 2i\}}.\\
\end{eqnarray*}
However 
\begin{eqnarray*}
\E |\Gamma_i^{-1/ \alpha(u)} - i^{-1/ \alpha(u)}|\mathbf{1}_{\{ \frac{i}{2} < \Gamma_i < 2i\}} & \leq & \frac{1}{i^{1/ \alpha(u)}} K_{u} \E |\frac{\Gamma_i}{i}-1|\\
& \leq & K_{u} \frac{1}{i^{1+\frac{1}{\alpha(u)}}} \\
\end{eqnarray*}
and
\begin{eqnarray*}
2^{j\eta}j^{2(1+ \eta)} \sum_{i=2^j}^{2^{j+1}-1}\limits \E |\Gamma_i^{-1/ \alpha(u)} - i^{-1/ \alpha(u)}|\mathbf{1}_{\{ \frac{i}{2} < \Gamma_i < 2i\}} & \leq & K j^{2(1+\eta)} 2^{j (\eta- \frac{1}{\alpha(u)})}.
\end{eqnarray*}
We thus obtain $\P \left( \liminf_j E_j\right) = 1$. As a consequence, $\sup_{v \in [0,1]}\limits \left|\sum_{i=1}^{+\infty}\limits \gamma_i (\Gamma_i^{-1/\alpha(u)}-i^{-1/ \alpha(u)} )\frac{\mathbf{1}_{[u,v]}(V_i)}{|v-u |^{\eta}} \right|< +\infty $ and finally
\begin{displaymath}
\sup_{v \in [0,1]}\limits \left| \frac{X(v,u)-X(u,u)}{|v-u |^{\eta}} \right|< +\infty 
\end{displaymath}
\Box
\begin{lem}\label{lemmajexpo}
For all $u \in (0,1)$, one has almost surely, for all $\eta \in (0, \frac{1}{\alpha(u)})$, 
\begin{displaymath}
\sup_{v \in [0,1]}  \frac{|X(v,u)-X(u,u)|}{|v-u |^{\eta}} < + \infty.
\end{displaymath}
\end{lem}

\noindent
{\bf Proof} 

\noindent
Fix $u \in (0,1)$. Lemma \ref{lemexpostab} yields that, for all $\eta \in (0, \frac{1}{\alpha(u)})$, we may choose an $\Omega_{\eta}$ having probability one and such that, on $\Omega_{\eta}$, $ \sup_{v \in [0,1]}\limits \left| \frac{X(v,u)-X(u,u)}{|v-u |^{\eta}} \right|< +\infty$. Thus, on $\Omega= \cap_{j \geq 0} \Omega_{\frac{1}{\alpha(u)}- \frac{1}{2^j}}$, which still has probability one, it holds that, for all $\eta \in (0,\frac{1}{\alpha(u)})$, $\sup_{v \in [0,1]}  \frac{|X(v,u)-X(u,u)|}{|v-u |^{\eta}} < + \infty$ \Box 

%
%

\bigskip

\noindent
{\bf Proof of Theorem \ref{expolevyexact}} 

\noindent
From Theorem \ref{expolevy}, we already know that $\mathcal{H}_u \leq \frac{1}{\alpha(u)}$.
To prove the reverse inequality, we treat separately the situations where $\alpha<1$ and $\alpha\geq 1$.
\begin{itemize}
\item Consider first the case $0< \alpha(u) <1$.
\end{itemize}
Write:
$$Y(v) - Y(u)  = X(v,v) - X(v,u) + X(v,u) - X(u,u).$$
By Lemma \ref{lemmajexpo}, we know that the Hölder regularity of $v \mapsto X(v,u) - X(u,u)$ at $u$ is almost
surely not smaller than $\frac{1}{\alpha(u)}$. Now, by applying the finite increments theorem to the functions
$t \mapsto C^{1/t}_{t}\Gamma_i^{-1/t}$, we get

\begin{eqnarray*}
	X(v,v) - X(v,u) & =& \sum_{i=1}^{\infty} \gamma_i \mathbf{1}_{[0,v]}(V_i) \left(C^{1/\alpha(v)}_{\alpha(v)}  \Gamma_i^{-1/\alpha(v)} - C^{1/\alpha(u)}_{\alpha(u)}  \Gamma_i^{-1/\alpha(u)}\right)\\
	& = & \left(\alpha(v)-\alpha(u)\right)\sum_{i=1}^{\infty} \gamma_i \mathbf{1}_{[0,v]}(V_i) 
	\left(CP(\alpha(w_i)) - C^{1/\alpha(w_i)}_{\alpha(w_i)}\frac{\log\Gamma_i}{\alpha(w_i)^2}\right) \Gamma_i^{-1/\alpha(w_i)}, 
\end{eqnarray*}
where, for each $i$, $w_i \in [u,v]$ (or $[v,u]$), and $CP$ denotes the derivative of the function 
$t \mapsto C^{1/t}_{t}$.
However,
\begin{eqnarray*}
|\sum_{i=1}^{\infty} \gamma_i \mathbf{1}_{[0,v]}(V_i) 
	\left(Cp(\alpha(w_i)) - \frac{\log\Gamma_i}{\alpha(w_i)^2}\right) \Gamma_i^{-1/\alpha(w_i)}| & \leq &
	\sum_{i=1}^{\infty} \left|CP(\alpha(w_i)) - C^{1/\alpha(w_i)}_{\alpha(w_i)}\frac{\log\Gamma_i}{\alpha(w_i)^2}\right| \Gamma_i^{-1/\alpha(w_i)} \\
	& \leq & K\sum_{i=1}^{\infty} \left(1+|\log\Gamma_i|\right) \left(\Gamma_i^{-1/c} + \Gamma_i^{-1/d}\right). 
\end{eqnarray*}
Thus the quantity $T(u,v) = \sum_{i=1}^{\infty} \gamma_i \mathbf{1}_{[0,v]}(V_i) 
	\left(CP(\alpha(w_i)) - C^{1/\alpha(w_i)}_{\alpha(w_i)}\frac{\log\Gamma_i}{\alpha(w_i)^2}\right) \Gamma_i^{-1/\alpha(w_i)}$ is, uniformly in $v$, almost
	surely finite and not 0. As a consequence, the function $v \mapsto X(v,v) - X(v,u) = (\alpha(u)-\alpha(v))T(u,v)$ has almost	surely the same Hölder exponent at $u$ as the function $v \mapsto \alpha(v)$ at $u$. If $\mathcal{H}_u^{\alpha} < \frac{1}{\alpha(u)}$, this entails that $Y$ has exponent
$\mathcal{H}_u^{\alpha}$ at $u$. If $\mathcal{H}_u^{\alpha} > \frac{1}{\alpha(u)}$, then the exponent of $Y$
at $u$ is at least $\frac{1}{\alpha(u)}$ and thus exactly $\frac{1}{\alpha(u)}$ by Theorem \ref{expolevy}.
 
\begin{itemize}
 \item Assume now that $1\leq \alpha(u) <2$. 
\end{itemize}
Let $\eta < \frac{1}{\alpha(u)}$ and $\delta \in (\eta, \frac{1}{\alpha(u)})$. Then:
\begin{displaymath}
 \frac{|Y(v)-Y(u)|}{|v-u|^{\eta}} \leq \frac{|X(v,v)-X(v,u)|}{|v-u|^{\eta}} + \frac{|X(v,u)-X(u,u)|}{|v-u|^{\eta}}.
\end{displaymath}
By lemma \ref{lemmajexpo}, there exists $K>0$ such that $\frac{|X(v,u)-X(u,u)|}{|v-u|^{\eta}} \leq K|v-u|^{\delta - \eta}$, and, by Lemma \ref{leminegtayl}, there exists $K>0$ such that $\frac{|X(v,v)-X(v,u)|}{|v-u|^{\eta}} \leq K|v-u|^{1 - \eta}.$
This entails $\lim_{v \rightarrow u} \frac{|Y(v)-Y(u)|}{|v-u|^{\eta}} = 0$ and
\begin{displaymath}
 \mathcal{H}_u \geq \frac{1}{\alpha(u)}
\end{displaymath}
\Box

\section{Assumptions}\label{Assum}

This section gathers the various conditions required on the considered processes so that
our results hold.
 
\begin{itemize}
    \item (C1) The family of functions $v \to f(t,v,x)$ is differentiable for all $(v,t)$ in $U^2$ and almost all $x$ in $E$. The derivatives of $f$ with respect to $v$ are denoted by $f'_v$.
    \item (C2) There exists $\delta > \frac{d}{c}-1$ such that :
\begin{displaymath}\label{kercond1f}
\sup_{t \in U}  \int_\bbbr \left[ \sup_{w \in U} (|f(t,w,x)|^{\alpha(w)}) \right]^{1+\delta} \hspace{0.1cm} \hat m(dx) < \infty.
\end{displaymath}
    \item (Cs2) There exists $\delta > \frac{d}{c} - 1$ such that :
\begin{displaymath}\label{kercond2sf}
\sup_{t \in U}  \int_\bbbr \left[ \sup_{w \in U} (|f(t,w,x)|^{\alpha(w)}) \right]^{1+\delta}  r(x)^{\delta} \hspace{0.1cm}  m(dx) < \infty.
\end{displaymath}
	\item (C3) There exists $\delta > \frac{d}{c}-1$ such that :
\begin{displaymath}
\sup_{t \in U}  \int_\bbbr \left[ \sup_{w \in U} (|f'_v(t,w,x)|^{\alpha(w)}) \right]^{1+\delta} \hspace{0.1cm} \hat m(dx) < \infty.
\end{displaymath} 
    \item (Cs3) There exists $\delta > \frac{d}{c} - 1$ such that :
\begin{displaymath}\label{kercond3sf}
\sup_{t \in U}  \int_\bbbr \left[ \sup_{w \in U} (|f'_v(t,w,x)|^{\alpha(w)}) \right]^{1+\delta} r(x)^{\delta} \hspace{0.1cm} m(dx) < \infty.
\end{displaymath}
    \item (Cs4) There exists $\delta > \frac{d}{c} - 1$ such that :
\begin{displaymath}\label{kercond5sf}
\sup_{t \in U}  \int_\bbbr \left[ \sup_{w \in U} \left[  \left|f(t,w,x) \log(r(x))  \right|^{\alpha(w)}  \right]\right]^{1+\delta} r(x)^{\delta} \hspace{0.1cm} m(dx) < \infty.
\end{displaymath}

	\item (C5) $X(t,u)$ (as a process in $t$) is localisable at $u$ with exponent $h(u) \in (h_-,h_+) \subset (0,1)$, with local form $X'_u(t,u)$, and  $u \mapsto h(u)$ is a $C^1$ function .
    \item (C6) There exists $K_U >0$ such that $\forall v \in U$, $\forall u \in U$, $\forall x \in \bbbr$,
    \begin{displaymath}
    \left| f(v,u,x) \right| \leq K_U.
    \end{displaymath}
    \item (C7) There exists $K_U >0$ such that $\forall v \in U$, $\forall u \in U$, $\forall x \in \bbbr$,
    \begin{displaymath}
    \left| f'_v(v,u,x) \right| \leq K_U.
    \end{displaymath}
     \item (C8) There exists a function $h$ defined on $U$, $\varep_0 \in ( 0,1)$ and $K_U >0$ such that $\forall r< \varep_0$, $\forall x \in \bbbr$,
    \begin{displaymath}
    \frac{1}{r^{h(t)-1/ \alpha(t)}}\left| f(t+r,t,x) -f(t,t,x) \right| \leq K_U.
    \end{displaymath}
     \item (Cu8) There exists a function $h$ defined on $U$ and $K_U >0$ such that $\forall v \in U$, $\forall u \in U$, $\forall x \in \bbbr$,
    \begin{displaymath}
    \frac{1}{|v-u|^{h(u)-1/ \alpha(u)}}\left| f(v,u,x) -f(u,u,x) \right| \leq K_U.
    \end{displaymath}
	\item (C9) There exists a function $h$ defined on $U$, $\varep_0 > 0$ and $K_U >0$ such that  $\forall r< \varep_0$,
\begin{displaymath}	
	\frac{1}{r^{h(t)\alpha(t)}} \int_{\bbbr} | f(t+r,t,x)-f(t,t,x)|^{\alpha(t)} m(dx)  \leq K_U.
\end{displaymath}
     \item (C10)  There exists a function $h$ defined on $U$ and $ p \in (\alpha(t),2)$, $p \geq 1$, such that for all $\varep > 0$, there exists $K_U >0$ such that,  $\forall r \leq \varep$,
	\begin{displaymath}
	 \frac{1}{r^{1+p(h(t)-\frac{1}{\alpha(t)})}} \int_{\bbbr} \left| f(t+r,t,x) - f(t,t,x) \right|^p m(dx)\leq K_U.
	\end{displaymath}
	\item (Cu10)  There exists a function $h$ defined on $U$, $ p \in (d,2)$, $p \geq 1$ and $ K_U >0$ such that $\forall v \in U$,$\forall u \in U$,
	\begin{displaymath}
	 \frac{1}{|v-u|^{1+p(h(u)-\frac{1}{\alpha(u)})}} \int_{\bbbr} \left| f(v,u,x) - f(u,u,x) \right|^p m(dx)\leq K_U.
	\end{displaymath}
	\item (C11) $\forall \varep > 0$, $\exists K_U >0$ such that, $\forall r \leq \varep$,
	\begin{displaymath}
	\int_{\bbbr} \left| f(t+r,t,x) \right|^2 m(dx) \leq K_U.
	\end{displaymath}
	\item (Cu11) There exists $ K_U >0$ such that $\forall  v \in U$, $\forall  u \in U$, 
	\begin{displaymath}
	\int_{\bbbr} \left| f(v,u,x) \right|^2 m(dx) \leq K_U.
	\end{displaymath}
    \item (C12) $\forall \varep > 0$, $\exists K_U >0$ such that $\forall r \leq \varep$,
	\begin{displaymath}
	\int_{\bbbr}  \left| f(t+r,t+r,x) \right|^2 m(dx) \leq K_U.
	\end{displaymath}
	
	\item (Cu12) There exists $ K_U >0$ such that $\forall  v \in U$, 
	\begin{displaymath}
	 \int_{\bbbr} \left| f(v,v,x) \right|^2 m(dx)  \leq K_U.
	\end{displaymath}
    \item (C13) 
    \begin{displaymath}
    \inf_{v \in U} \int_{\bbbr} f(v,v,x)^2 m(dx) >0.
    \end{displaymath}
	\item (C14) There exists a function $h$ and a positive function $g$ defined on $U$ such that 
	\begin{displaymath}
	\lim_{r \rightarrow 0} \frac{1}{r^{1+2(h(t)-1/\alpha(t)})} \int_{\bbbr} \left( f(t+r,t,x) - f(t,t,x) \right)^2  m(dx) = g(t).
	\end{displaymath}
	\item (Cu14) There exists a function $h$ and a positive function $g$ defined on $U$ such that 
    \begin{displaymath}
\lim_{r \rightarrow 0} \sup_{t \in U} \left| \frac{1}{r^{1+2(h(t)-1/\alpha(t)})} \int_{\bbbr} \left( f(t+r,t,x) - f(t,t,x) \right)^2  m(dx) - g(t) \right|  = 0.
\end{displaymath}

	\item (C15) $\forall \varep > 0$, $\exists K_U >0$ such that $\forall r \leq \varep$,
	\begin{displaymath}
	 \frac{1}{|r|^2} \int_{\bbbr} \left| f(t+r,t+r,x) - f(t+r,t,x) \right|^2 m(dx)\leq K_U.
	\end{displaymath}
	
	\item (Cu15) $\exists K_U >0$ such that, $\forall v \in U$, $\forall u \in U$, 
	\begin{displaymath}
	 \frac{1}{|v-u|^2} \int_{\bbbr} \left| f(v,v,x) - f(v,u,x) \right|^2 m(dx)\leq K_U.
	\end{displaymath}

\end{itemize}


\begin{thebibliography}{abc-17}

\bibitem{AL2}  {\sc Ayache, A. and L\'{e}vy V\'{e}hel, J.} (2000).
The generalized
multifractional Brownian motion. \emph{Stat. Inference 
Stoch. Process.} {\bf 3,} 7--18.

\bibitem{BE} {\sc Von Bahr, B. and Essen, C.G.} (1965).
Inequalities for the rth Absolute Moment of a Sum of Random Variables, 1 <=r <= 2
\emph{The Annals of Mathematical Statistics} {\bf 36,} (1), 299--303.

\bibitem{BJR}{\sc  Benassi, A., Jaffard, S. and Roux, D.} (1997). Gaussian
processes and pseudodifferential elliptic operators.
\emph{Rev. Mat. Iberoamericana } {\bf 13,} 19--89.

\bibitem{BJP} {\sc Bentkus, V., Juozulynas, A. and Paulauskas, V.} (2001).
L\'{e}vy-LePage series representation of stable vectors: convergence in variation.
\emph{J. Theo. Prob.}(2) {\bf 14,} (4) 949--978.

\bibitem{Fal5} {\sc Falconer, K.J.} (2002). Tangent fields and the local structure
of random fields. \emph{J. Theoret. Probab.} {\bf 15,}
731--750.

\bibitem{Fal6} {\sc Falconer, K.J.} (2003). The local structure of random processes.
 \emph{J. London Math. Soc.}(2) {\bf 67,}
657--672. 

\bibitem{FLGLV} {\sc Falconer, K.J., Le Gu\'{e}vel, R. and L\'{e}vy V\'{e}hel, J.} (2009).
Localisable moving average stable and multistable processes, {\it Stochastic Models}, 2009, no.4, 648--672.

\bibitem{FLV} {\sc Falconer, K.J. and L\'{e}vy V\'{e}hel, J.} (2008).
{\it Multifractional, multistable, and other processes
with prescribed local form}, {\it J. Theoret. Probab.},
DOI 10.1007/s10959-008-0147-9.

\bibitem{FL} {\sc Falconer, K.J., and Lining, L.} (2009).
Multistable random measures and multistable processes, preprint.

\bibitem{FK} {\sc Ferguson, T.S. and Klass, M.J.} (1972).
A representation of independent increment processes without Gaussian components.
\emph{Ann. Math. Stat.} {\bf 43,} 1634--1643.

\bibitem{EH}  {\sc Herbin, E.} (2006). From N-parameter fractional Brownian motions
to N-parameter multifractional Brownian motion. \emph{Rocky
Mountain J. Math.} {\bf 36,} 1249--1284.

\bibitem{EHJLV}{\sc Herbin, E. and Lévy Véhel, J.} (2009). Stochastic 2 micro-local analysis.
\emph{Stoch. Proc. Applic.} {\bf 119}(7), 2277--2311.

\bibitem{K40} {\sc Kolmogorov, A.N.} (1940). Wienersche Spiralen und einige
andere interessante Kurven in Hilbertchen Raume, \emph{Doklady},
{\bf 26}, 115--118.

\bibitem{LGLV} {\sc Le Gu\'{e}vel, R. and L\'{e}vy V\'{e}hel, J.} (2009).
{\it A Ferguson - Klass - LePage series representation of
multistable multifractional motions and related processes}, preprint.
Available at http://arxiv.org/abs/0906.5042.

\bibitem{LP1} {\sc Le Page, R.} (1980).
Multidimensional infinitely divisible variables and processes. I. Stable case
{\emph Tech. Rep. 292, Dept. Stat}, Stanford Univ.

\bibitem{LP2}{\sc Le Page, R.} (1980).
Multidimensional infinitely divisible variables and processes. II
{\emph Probability in Banach Spaces III}
Lecture notes in Math. {\bf 860} 279--284, Springer, New York.

\bibitem{LTal} {\sc Ledoux, M. and Talagrand, M.} (1996).
Probability in Banach spaces. {\it Springer-Verlag.}

\bibitem{MV} {\sc Mandelbrot, B.B. and Van Ness, J.} (1968). Fractional Brownian
motion, fractional noises and applications.
 \emph{SIAM Rev.} {\bf 10,} 422--437.

\bibitem{PL} {\sc  Peltier, R.F. and L\'{e}vy V\'{e}hel, J.} (1995).
 Multifractional Brownian motion: definition and preliminary results.
 \emph{Rapport de recherche de l'INRIA}, No. 2645. Available
at:
http://www-rocq1.inria.fr/fractales/index.php?page=publications

\bibitem{Pet} {\sc Petrov, V.} (1995). {\it Limit Theorems of Probability Theory},
Oxford Science Publication.

\bibitem{JR}{\sc Rosinski, J.} (1990).
On Series Representations of Infinitely Divisible Random Vectors
\emph{Ann. Probab.}, {\bf 18,} (1) 405--430.

\bibitem{ST} {\sc Samorodnitsky, G. and Taqqu, M.S.} (1994).
{\it Stable Non-Gaussian Random Processes},
Chapman and Hall.

\bibitem{ST1} {\sc Stoev, S. and Taqqu, M.S.} (2004).
Stochastic properties of the linear multifractional stable motion.
\emph{Adv.  Appl. Probab.}, {\bf 36,}
1085--1115

\bibitem{ST2} {\sc Stoev, S. and Taqqu, M.S.} (2005).
Path properties of the linear multifractional stable motion. \emph{
Fractals}, {\bf 13,} 157--178.



\end{thebibliography}
\end{document}